\documentclass[12pt,twoside]{article}
\usepackage{a4wide,color,bm,amssymb,amsmath,hyperref}
\newtheorem{theorem}{Theorem}
\newtheorem{proposition}{Proposition}

\newtheorem{definition}{Definition}

\newtheorem{lemma}{Lemma}

\renewcommand{\Bbb}[1]{\mathbb{#1}}
\newcommand{\C}{{\Bbb C}}         
\newcommand{\N}{{\Bbb N}}         
\newcommand{\Q}{{\Bbb Q}}         
\newcommand{\R}{{\Bbb R}}         
\newcommand{\Rp}{{\Bbb R}^{+}}    
\newcommand{\Z}{{\Bbb Z}}         
\newcommand{\Rmn}{\R^{m\times n}}
\newcommand{\xmn}{\bm X}

\newcommand{\txmn}{{}^t\!\xmn}
\newcommand{\wx}{w(\xmn)}
\newcommand{\wxth}{w(\xmn,\bm\theta)}
\newcommand{\hwx}{\widehat w(\xmn)}
\newcommand{\hwxth}{\widehat w(\xmn,\bm\theta)}
\newcommand{\ws}{w^\times}
\newcommand{\wsx}{\ws(\xmn)}
\newcommand{\wsxth}{\ws(\xmn,\bm\theta)}
\newcommand{\bth}{\bm\theta}
\newcommand{\mn}{m\times n}

\newcommand{\De}{\Delta}

\newcommand{\La}{\Lambda}

\newcommand{\p}{\psi}
\renewcommand{\hom}{{\rm H}}
\newcommand{\inh}{{\rm I}}

\newcommand{\cA}{{\cal A}}

\newcommand{\cC}{{\cal C}}
\newcommand{\cD}{{\cal D}}
\newcommand{\cE}{{\cal E}}

\newcommand{\cG}{{\cal G}}

\newcommand{\cL}{{\cal L}}
\newcommand{\cM}{{\cal M}}
\newcommand{\cN}{{\cal N}}

\newcommand{\cP}{{\cal P}}

\newcommand{\cV}{{\cal V}}

\newcommand{\Bn}{B_{\mathrm{N}}}

\newcommand{\ttt}{{\vv t}}
\newcommand{\TTT}{{\vv T}}
\newcommand{\AAA}{\mathcal{A}}
\newcommand{\Cta}{{\cC}_{\ttt,\alpha}}

\newcommand{\Gta}{\cG_{\ttt,\alpha}}
\newcommand{\Dta}{\cD_{\ttt,\alpha}}
\newcommand{\Nta}{\cN_{\ttt,\alpha}}
\newcommand{\iDta}{\inh_\ttt(\alpha,\p_\ttt)}

\newcommand{\Lad}{\La_{\mathrm{D}}}
\newcommand{\Lan}{\La_{\mathrm{N}}}

\newcommand{\PI}{{\textstyle\prod}}
\newcommand{\lan}{\langle}
\newcommand{\ran}{\rangle}

\newcommand{\ve}{\varepsilon}
\newcommand{\supp}{\operatorname{supp}}
\newcommand{\Supp}{\bm{S}}
\newcommand{\diam}{\operatorname{diam}}
\newcommand{\dist}{\operatorname{dist}}

\newcommand{\vv}[1]{{\mathbf{#1}}}

\renewcommand{\le}{\leqslant}
\renewcommand{\ge}{\geqslant}

\newcommand{\bnz}{\smallsetminus\{\vv0\}}
\newcommand{\s}[1]{\sigma(\vv#1)}

\newcommand{\qqand}{\qquad\text{and}\qquad}
\newcommand{\wh}{\widehat}

\begin{document}

\title{ An Inhomogeneous Transference Principle  \\ and Diophantine Approximation}

\author{Victor Beresnevich\footnote{EPSRC Advanced Research Fellow, grant EP/C54076X/1}
\\ {\small\sc York} \and Sanju Velani\footnote{Research supported by EPSRC  grant EP/E061613/1 and
INTAS grant 03-51-5070}
\\ {\small\sc York}}

\bigskip

\date{{\small\it Dedicated to Vasili  Bernik -- spasibo Basil! }}

\maketitle


\begin{abstract}
In a landmark paper \cite{Kleinbock-Margulis-98:MR1652916},
D.Y.\;Kleinbock and G.A.\;Margulis  established the fundamental
Baker-Sprindzuk conjecture on homogeneous Diophantine approximation
on manifolds. Subsequently, there has been dramatic progress in this
area of research. However, the techniques developed to date  do not
seem to be applicable to \emph{inhomogeneous}\/ approximation.
Consequently, the theory of inhomogeneous Diophantine approximation
on manifolds remains essentially non-existent.

In this paper we develop an approach that enables us to transfer
homogeneous statements to inhomogeneous ones. This is rather
surprising as the inhomogeneous theory  contains the homogeneous
theory and so is more general.  As a consequence, we establish the
inhomogeneous analogue of the Baker-Sprindzuk conjecture.
Furthermore, we prove a complete inhomogeneous version of the
profound theorem of Kleinbock, Lindenstrauss $\&$  Weiss
\cite{Kleinbock-Lindenstrauss-Weiss-04:MR2134453} on the
extremality of friendly measures. The results obtained in this
paper constitute the first step towards developing a coherent
inhomogeneous theory for manifolds in line with the homogeneous
theory.
\end{abstract}

\newpage


\section{Introduction}

The metrical theory of Diophantine approximation on manifolds dates
back to the nineteen thirties with a conjecture of K.~Mahler
\cite{Mahler-1932b} in transcendence theory. The conjecture was
easily seen to be equivalent to a metrical Diophantine approximation
problem restricted to the Veronese curves ${\cal V}_n := \{
(x,\dots,x^n)$, $x\in\R \}$. Mahler's conjecture remained a key open
problem in metric number theory for over 30 years and was eventually
solved  by Sprindzuk \cite{Sprindzuk-1969-Mahler-problem}. Moreover,
its solution led Sprindzuk \cite{Sprindzuk-1980-Achievements} to
make an important general conjecture which we shall  shortly
describe. The conjecture has  been established  by Kleinbock $\&$
Margulis in their landmark paper
\cite{Kleinbock-Margulis-98:MR1652916}. The main result of this
paper establishes a complete inhomogeneous version of the theorem of
Kleinbock $\&$ Margulis and indeed its generalisation to friendly
measures \cite{Kleinbock-Lindenstrauss-Weiss-04:MR2134453}. In order
to describe these fundamental conjectures and results it is
convenient to introduce the notion of Diophantine exponents.

\subsection{Exponents of Diophantine approximation}

Let $m,n\in\N$ and $\Rmn$ be the set of all $\mn$ real matrices.
Given $\xmn\in\Rmn$ and $\bth\in\R^{m}$, let $\wxth$ be the supremum
of $w\ge0$ such that for arbitrarily large $Q>1$ there is a $\vv
q=(q_1,\dots,q_n)\in\Z^n\bnz$ satisfying
\begin{equation}\label{e:001}
 \|\xmn\vv q+\bth\|^m<Q^{-w}\qqand|\vv q|^n\le Q\,,
\end{equation}
where $|\vv q|:= \max\{|q_1|, \ldots, |q_n|\} $ is the supremum norm
and $\|\cdot\|$ is the distance to the nearest integer point. Here
and elsewhere $\vv q\in\Z^n$ and $\bth\in\R^m$ are treated as
columns. It follows that whenever $\wxth$ is finite, the inequality
$$
\|\xmn\vv q+\bth\|^m<|\vv q|^{-wn}
$$
has infinitely many solutions $\vv q\in\Z^n$ if $w<\wxth$ and has
at most finitely many solutions $\vv q\in\Z^n$ if $w>\wxth$.
Further, let $\wsxth$ be the supremum of $w\ge0$ such that for
arbitrarily large $Q>1$ there is a $\vv
q=(q_1,\dots,q_n)\in\Z^n\bnz$ satisfying
\begin{equation}\label{e:002}
 \PI\lan\xmn\vv q+\bth\ran<Q^{-w}\qqand\PI_+(\vv q)\le Q\,,
\end{equation}
where
$$
\text{$\PI\vv y:=\PI(\vv y)=\prod_{j=1}^m|y_j|$\qquad and\qquad
$\PI_+(\vv q):=\prod_{i=1}^n\max\{1,|q_i|\} $}
$$
for $\vv y=(y_1,\dots,y_m)$. Also $\lan\vv y\ran$ denotes the unique
point in $[-1/2,1/2)^m$ congruent to $\vv y\in\R^m$ modulo $\Z^m$.
Thus, $\|\cdot\|=|\lan\ \cdot\ \ran|$. It follows  that whenever
$\wsxth$ is finite, the inequality
$$
\PI\lan\xmn\vv q+\bth\ran<\PI_+(\vv q)^{-w}
$$
has infinitely many solutions $\vv q\in\Z^n$ if $w<\wsxth$ and has
at most finitely many solutions $\vv q\in\Z^n$ if $w>\wsxth$.

The homogeneous theory of Diophantine approximation corresponds to
the special case of $\bth=\vv0$ in the above inequalities. In this
case the associated  homogeneous exponents will be denoted by
$\wx:=w(\xmn,\vv0)$ and $\wsx:=w^*(\xmn,\vv0)$. A trivial
consequence of Dirichlet's theorem \cite{Schmidt-1980}, or simply
the `pigeon-hole principle', is that
\begin{equation}\label{e:003}
 \wx\ge 1\qquad\text{for all }\xmn\in\Rmn\,.
\end{equation}
Also it is readily seen that (\ref{e:001}) implies (\ref{e:002}) and
therefore
\begin{equation}\label{e:004}
 \wsxth\ge \wxth\qquad\text{for all }\xmn\in\Rmn\text{ and all }\bth\in\R^m\,.
\end{equation}

The Diophantine exponents can in principle be infinite.
Nevertheless, a relatively straightforward consequence of the
Borel-Cantelli lemma from probability theory is that (\ref{e:003})
is reversed for almost all $\xmn\in\Rmn$ with respect to Lebesgue
measure on $\Rmn$ and moreover that
\begin{equation}\label{e:005}
 \wsx=1\qquad\text{for almost all }\xmn\in\Rmn\,.
\end{equation}

\noindent For completeness, we mention that $\xmn\in\Rmn$ is said
to be very well approximable  (see
\cite{Kleinbock-Margulis-98:MR1652916,Schmidt-1980})
 if $ \wx>1 $ and multiplicatively
very well approximable  (see \cite{Kleinbock-Margulis-98:MR1652916})
if $ \wsx>1 $.  Note that in view of the discussion above, the
corresponding sets of very well approximable and multiplicatively
very well approximable points are of zero Lebesgue measure on
$\Rmn$.

\subsection{Homogeneous theory}

Sprindzuk \cite{Sprindzuk-1980-Achievements} conjectured that
(\ref{e:005}) remains true when $\xmn$ is restricted to any analytic
non-degenerate submanifold $\cM$ of $\R^n$ identified with either
columns $\R^{n\times1}$ (\emph{simultaneous Diophantine
approximation}) or rows $\R^{1\times n}$ (\emph{dual Diophantine
approximation}) with respect to the Riemannian measure on $\cM$.
This conjecture had been previously stated by A.\,Baker
\cite{Baker-75:MR0422171} for Veronese curves ${\cal V}_n := \{
(x,\dots,x^n)$, $x\in\R \}$. Essentially, non-degenerate manifolds
are smooth submanifolds of $\R^n$ which are sufficiently curved so
that they deviate from any hyperplane with a `power law'
\cite{Beresnevich-02:MR1905790}. For the formal definition see
\cite{Kleinbock-Margulis-98:MR1652916}. Any real, connected analytic
manifold not contained in a hyperplane of $\R^n$ is non-degenerate.
For a planar curve, the non-degeneracy condition is simply
equivalent to the condition that the curvature is non-vanishing
almost everywhere.

\medskip

\noindent\textbf{Baker-Sprindzuk conjecture.} {\it For any analytic
non-degenerate submanifold $\cM$ of\/ $\R^n$
\begin{equation}\label{e:006}
 \wsx=1\qquad\text{for almost all }\xmn\in
 \cM\,.
\end{equation}
}

Note that in view of (\ref{e:003}) and (\ref{e:004}), the
Baker-Sprindzuk conjecture implies that
\begin{equation}\label{e:007}
 w(\xmn)=1\qquad\text{for almost all }\xmn\in \cM\,.
\end{equation}

\noindent In fact, this weaker statement  also  appears as a
formal conjecture in \cite{Sprindzuk-1980-Achievements}.  In the
case of the Veronese curves ${\cal V}_n$, (\ref{e:007}) reduces to
Mahler's problem \cite{Mahler-1932b} and statement (\ref{e:006})
reduces to the specific conjecture of A.~Baker mentioned above.
Manifolds that satisfy (\ref{e:007}) are referred to as
\emph{extremal} and manifolds that satisfy (\ref{e:006}) are
referred to as \emph{strongly extremal}. To be precise, either
notion of extremality actually takes on two forms depending on
whether $\R^n$ is identified with $\R^{n\times 1}$ or $\R^{1\times
n}$.
  However, by
Khintchine's transference principle both forms are equivalent and it
is pointless to distinguish between them. A priori, this is not the
case when considering `extremality' within the inhomogeneous setting
-- see \S\ref{inhomtheory}.

Until 1998, progress on the Baker-Sprindzuk conjecture was limited
to special classes of manifolds -- see \cite{BernikDodson-1999,
Sprindzuk-1979-Metrical-theory}. Restricting our attention to
non-degenerate curves, the weaker form of the conjecture
corresponding to  (\ref{e:007}) had been established for planar
curves by W.M.~Schmidt \cite{Schmidt-1964b}
 and by Beresnevich $\&$ Bernik \cite{Beresnevich-Bernik-96:MR1387861} for curves in
 $\R^3$.  The actual conjecture  had  only been established  for
the Veronese curves ${\cal V}_n$ with $n\le 3$ by Yu
\cite{Yu-81:MR620731} and with $n=4$ by Bernik \& Borbat
\cite{Bernik-Borbat-97:MR1639563}. In their ground breaking work,
Kleinbock $\&$ Margulis \cite{Kleinbock-Margulis-98:MR1652916}
established the Baker-Sprindzuk conjecture in full generality and
moreover removed the `analytic' assumption.
\medskip

\noindent\textbf{Theorem KM } {\it Any non-degenerate submanifold of
$\R^n$ is strongly extremal.}

\medskip

 The work of Kleinbock $\&$  Margulis has led to various
generalisations of the Baker-Sprindzuk conjecture. Kleinbock
\cite{Kleinbock-03:MR1982150} has established that non-degenerate
submanifolds of strongly extremal affine subspaces of $\R^n$ are
strongly extremal. He has also shown that non-degenerate complex
analytic manifolds are strongly extremal
\cite{Kleinbock-04:MR2094125}. Kleinbock $\&$  Tomanov
\cite{Kleinbock-Tomanov-07:MR2314053} have generalised Theorem KM
to the $S$-arithmetic setting. Kleinbock, Lindenstrauss $\&$ Weiss
\cite{Kleinbock-Lindenstrauss-Weiss-04:MR2134453} have
revolutionised the notion of extremality by introducing the
concept of measures being extremal rather than sets.  Let $\mu$ be
a measure supported on a subset of $\Rmn$. We say that $\mu$ is
\emph{extremal} if $\wx=1$ for $\mu$-almost every point
$\xmn\in\Rmn$.  In other words,  the set of $\xmn\in\Rmn$ for
which $\wx>1$ is of $\mu$-measure zero. We say that $\mu$ is
\emph{strongly extremal} if $\wsx=1$ for $\mu$-almost every point
$\xmn\in\Rmn$. Furthermore, if $\mu$ is a measure on $\R^n$ then
$\mu$ is \emph{(strongly) extremal}\/ if it is (strongly) extremal
through the identification of $\R^n$ with either $\R^{n\times1}$
or $\R^{1\times n}$. In view of Khintchine's transference
principle, there is no difference which representation of $\R^n$
is taken.

The following constitutes  the main result of Kleinbock,
Lindenstrauss $\&$ Weiss
\cite{Kleinbock-Lindenstrauss-Weiss-04:MR2134453}.

\medskip

\noindent\textbf{Theorem KLW } {\it Any friendly measure on\/ $\R^n$
is strongly extremal.}

\medskip

\noindent The definition of {\em friendly }  measures is given in
\S\ref{frm}. At this point it suffices to say that friendly
measures form a large and natural class of measures on $\R^n$
including Riemannian measures supported on non-degenerate
manifolds, fractal measures supported on self-similar  sets
satisfying the open set condition (e.g. regular Cantor sets, Koch
snowflake, Sierpinski gasket)  and conformal measures supported on
limit sets of Kleinian groups.   In view of the former we have
that
$$
\text{Theorem KLW \quad$\Longrightarrow$\quad Theorem KM. }
$$

\subsection{Inhomogeneous theory \label{inhomtheory}}

The central goal of this paper is to establish the inhomogeneous
analogue of the Baker-Sprindzuk conjecture. Naturally, we
 begin by introducing the notion of extremality in the
inhomogeneous theory of Diophantine approximation.

\begin{definition}[Measures on $\Rmn$]\rm Let $\mu$ be a
measure supported on a subset of $\Rmn$. We say that $\mu$ is
\emph{inhomogeneously extremal} if for all $\bth\in\R^{m}$
\begin{equation}\label{e:008}
\wxth=1\qquad\text{for $\mu$-almost all $\xmn\in\Rmn$.}
\end{equation}
We say that $\mu$ is \emph{inhomogeneously strongly extremal} if for
all $\bth\in\R^{m}$
\begin{equation}\label{e:009}
    \wsxth=1\qquad\text{for $\mu$-almost all $\xmn\in\Rmn$.}
\end{equation}
\end{definition}

\noindent These notions of extremality naturally generalise the
homogenous ones which only require (\ref{e:008}) and (\ref{e:009})
to hold for $\bth=\vv0$. A remark regarding the use of the word
\emph{strongly} in the definition of `inhomogeneously strongly
extremal' is in order. In the homogeneous case, (\ref{e:003}) and
$(\ref{e:004})$ with $\bth=\vv0 $  show that strong extremality
implies extremality -- exactly as one would expect.  In the
inhomogeneous case there is no analogue of (\ref{e:003}) and it is
not at all obvious  that strong extremality implies extremality.
However, the following result established in \S\ref{sec:1B}
justifies the use of the word `strongly' even in the inhomogeneous
case.

\begin{proposition}\label{prop1}
Let $\mu$ be a measure on $\Rmn$. Then

\medskip

 \centerline{$\mu$ is inhomogeneously strongly extremal $\quad\Longrightarrow\quad$ $\mu$ is inhomogeneously
extremal\,.}
\end{proposition}

As already mentioned, there are two different forms of Diophantine
approximation when approximating points in $\R^n$ depending on
whether $\R^n$ is identified with $\R^{n\times 1}$ or $\R^{1\times
n}$. The identification with the former corresponds to the
simultaneous form  and the latter corresponds to the dual form. As a
consequence of Khintchine's transference principle, the two forms of
approximation  lead to equivalent notions of extremality in the
homogeneous case. However, Khintchine's transference principle is
not applicable in the inhomogeneous case and the simultaneous and
dual forms of extremality are not necessarily equivalent.
Consequently, the two forms of extremality need to be considered
separately.

\begin{definition}[Measures on $\R^n$]\label{def2}\rm
Let $\mu$ be a measure supported on a subset of $\R^n$. If $\mu$ is
inhomogeneously (strongly) extremal on $\R^{1\times n}$ we say that
$\mu$ is \emph{dually inhomogeneously (strongly) extremal}. If $\mu$
is inhomogeneously (strongly) extremal on $\R^{n\times 1}$ we say
that $\mu$ is \emph{simultaneously inhomogeneously (strongly)
extremal}. If $\mu$ is both dually and simultaneously
inhomogeneously (strongly) extremal then we simply say that $\mu$ is
\emph{inhomogeneously (strongly) extremal}.
\end{definition}

Naturally, a manifold $\cM\subset\R^n$ is called
\emph{inhomogeneously (strongly) extremal}\/ if the Riemannian
measure on $\cM$ is inhomogeneously (strongly) extremal. The
 Veronese
curves $\cV_n$ have been shown to be dually  inhomogeneously
extremal in the real \cite{Bernik-Dickinson-Dodson-98:MR1712779},
complex \cite{Ustinov-06:MR2229813}, $p$-adic
\cite{Bernik-Dickinson-Yuan-1999-inhom, Ustinov-05:MR2210400} and
`mixed' \cite{Bernik-Kovalevskaya-06:MR2298911} cases. These results
are natural generalisations of Mahler's conjecture to the
inhomogeneous setting. Most recently, Badziahin \cite{Bodyagin} has
extended Schmidt's homogeneous result \cite{Schmidt-1964b} by
showing that non-degenerate planar curves are dually inhomogeneously
extremal. Strikingly this constitutes the only known result beyond
the Veronese curves. However, one would expect that an inhomogeneous
analogue of Theorem~KM holds in full generality.

\medskip

\noindent\textbf{Inhomogeneous Baker-Sprindzuk conjecture.} {\it
Any non-degenerate submanifold  of\/ $\R^n$ is inhomogeneously
strongly extremal. }

\medskip

\noindent In view of Proposition~\ref{prop1}, this implies the
weaker conjecture that {\it any non-degenerate submanifold $\cM$ of
$\R^n$ is inhomogeneously extremal}. As discussed above the weaker
conjecture is known to be true for Veronese curves. Regarding the
stronger conjecture nothing is known. Nevertheless, given
Theorem~KLW, it is natural to broaden the conjecture to
friendly measures.

\medskip

\noindent\textbf{Conjecture.} \emph{Any friendly measure on $\R^n$
is inhomogeneously strongly extremal.}

\subsection{Statement of results}\label{svresults}

Let $\mu$ be a non-atomic, locally finite, Borel measure on $\Rmn$.
Obviously, we have that
$$
\text{$\mu$ is inhomogeneously (strongly) extremal
\quad$\Longrightarrow$\quad$\mu$ is (strongly) extremal.}
$$
The right hand side corresponds to the special choice of $\bth=\vv0$
in the definition of inhomogeneously (strongly) extremal. In this
paper we show that the above implication is reversed for a large
class of measures.

\begin{theorem}\label{t1}
Let $\mu$ be a measure on $\Rmn$.
\begin{itemize}
\item[{\rm (A)}] If $\mu$ is contracting almost everywhere then

\smallskip

 \centerline{$\mu$ is extremal $\iff$ $\mu$ is inhomogeneously extremal.}

\item[{\rm(B)}] If $\mu$ is strongly contracting almost everywhere then

\smallskip

 \centerline{$\mu$ is strongly extremal $\iff$ $\mu$ is inhomogeneously strongly extremal.}
\end{itemize}
\end{theorem}

 So as to avoid introducing various
technical notions at this point, the definition of contracting
measures is postponed to the next section. It suffices to say that
friendly measures on $\R^n$ fall within the class of strongly
contracting measures. Thus, Theorem \ref{t1}B (part (B) of Theorem
\ref{t1}) together with Theorem KLW establishes the conjecture
stated above for friendly measures.

\begin{theorem}\label{t2}
Any friendly measure on $\R^n$ is inhomogeneously strongly extremal.
\end{theorem}

 Riemannian measures supported on  non-degenerate manifolds are
know to be friendly
\cite{Kleinbock-Lindenstrauss-Weiss-04:MR2134453}.  Thus,
Theorem~\ref{t2} gives a complete inhomogeneous analogue of
Theorem~KM and thereby settles the inhomogeneous Baker-Sprindzuk
conjecture.

\begin{theorem}\label{t3}
Any non-degenerate submanifold of\/ $\R^n$ is inhomogeneously
strongly extremal.
\end{theorem}

\bigskip

It is worth mentioning that the class of contracting measures is
not limited to friendly measures.  To illustrate this we restrict
our attention to simultaneous Diophantine approximation. Thus,
$\R^n$ is identified with $\R^{n\times1}$. In \S\ref{svdiffman},
we show that the Riemannian measure on an arbitrary differentiable
submanifold $\cM$ of $\R^n$  falls within the class of contracting
measures and indeed within the class of strongly contracting
measure if a mild condition is imposed on $\cM$. However, affine
subspaces of $\R^n$ are differentiable manifolds  and for obvious
reasons they do not support friendly measures. Specializing
Theorem \ref{t1} to differentiable manifolds gives the following
statement.

\begin{theorem}\label{t4}
Let $\cM$ be a differentiable submanifold of $\R^{n}$. Then\\[0.9ex]
\phantom{}{\rm(A)}~~$\cM$ is extremal $\iff$ $\cM$ is simultaneously
inhomogeneously
 extremal\,.\\[2ex]
\phantom{}Furthermore, suppose that at almost every point on $\cM$
the tangent plane is not orthogonal to
any of the coordinate axes. Then\\[0.9ex]
\phantom{}{\rm(B)}~~$\cM$ is strongly extremal $\iff$ $\cM$ is
simultaneously inhomogeneously strongly extremal.
\end{theorem}

\noindent Recall, that the left hand side of the above
implications are homogeneous statements and the notions of
simultaneously and dually (strongly) extremal coincide.  Examples
of (strongly) extremal differentiable submanifolds that do not
fall within the remit of Theorem  \ref{t2}  are given in
\cite{Beresnevich-Bernik-Dickinson-Dodson-00:MR1820985,Kleinbock-03:MR1982150,Schmidt-1964b}.
Thus, Theorem~\ref{t4} is not vacuous.

\medskip

The following diagram summarizes the connections between the various
notions of extremality for strongly contracting measures on $\Rmn$.

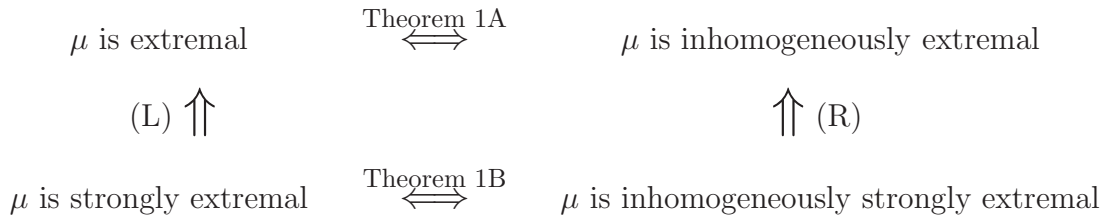
$$
\begin{array}{ccc}
\text{$\mu$ is extremal} \hspace{2ex} &
\text{\large$\stackrel{\text{Theorem \ref{t1}A}}{\iff}$} &
\hspace{2ex} \text{$\mu$ is inhomogeneously extremal} \\[2ex]
\text{(L) \LARGE$\Uparrow$} &  & \text{{\LARGE$\Uparrow$} (R)}\\[2ex]
\text{$\mu$ is strongly extremal}  \hspace{2ex} &
\text{\large$\stackrel{\text{Theorem \ref{t1}B}}{\iff}$} &
\hspace{2ex} \text{$\mu$ is inhomogeneously strongly extremal}
\end{array}
$$

\medskip

\noindent As mentioned in \S\ref{inhomtheory} the implication (L) is
well known. Thus, for contracting measures the implication (R)
follows via Theorem \ref{t1} and is independent of
Proposition~\ref{prop1}.

\bigskip

The upshot of Theorem \ref{t1} is that it enables us to transfer
homogeneous extremality statements to inhomogeneous ones. This
`inhomogeneous transference' is rather surprising as the
inhomogeneous theory contains the homogeneous theory and so is more
general. In \S\ref{HIT}, we develop an abstract
framework within which we establish a general inhomogeneous
transference principle -- namely Theorem \ref{t5}. The key step
in establishing Theorem \ref{t1} follows as an application
of this inhomogeneous transference principle.

\medskip

\noindent{\em Remark. } A direct and self-contained proof of
Theorem \ref{t4}A can be found in our recent article
\cite{Beresnevich-Velani-Moscow}. The main motivation behind
\cite{Beresnevich-Velani-Moscow} is to foreground and
significantly simplify the key ideas involved in establishing the
inhomogeneous transference principle of \S\ref{HIT}.  Indeed,
anyone interested in the proof of Theorem \ref{t5} and thus
Theorem \ref{t1} may find it useful first to look at
\cite{Beresnevich-Velani-Moscow}.

\section{Contracting and friendly measures}\label{frm}

In this section we start by formally introducing the class of
`contracting' measures alluded to in Theorem \ref{t1} above. We then
show that friendly measures are contracting and that the Riemannian
measure on a differentiable submanifold of $\R^{n\times1}$ is
contracting. This establishes Theorems~\ref{t2} and \ref{t4} from
Theorem~\ref{t1}.

\subsection{Contracting measures}\label{contr}

We begin by recalling some standard notions. If $B$ is a ball in a
metric space $\Omega$ then $cB$ denotes the ball with the same
centre as $B$ and radius  $c$ times the radius of $B$. A measure
$\mu$ on $\Omega$ is \emph{non-atomic} if the measure of any point
in $\Omega$ is zero. The \emph{support} of $\mu$ is the smallest
closed set $\Supp$ such $\mu(\Omega\setminus \Supp) = 0$. Also,
recall that $\mu$ is \emph{doubling} if there is a constant
$\lambda>1$ such that for any ball $B$ with centre in $\Supp$
\begin{equation}\label{e:010}
  \mu\big(2B\big)\ \le\ \lambda\ \mu\big(B\big)\,.
 \end{equation}

\noindent The class of contracting measures $\mu$ is defined via the
behavior of $\mu$ near planes in $\Rmn$ .  More precisely, the
planes are given by
\begin{equation}\label{e:011}
\cL_{\vv a,\vv b}:=\{\xmn\in\Rmn:\xmn\vv a+\vv b=\vv0\}
\qquad\text{with}\quad\vv a\in\R^n,\quad|\vv
a|_{2}=1\quad\text{and}\quad\vv b\in\R^m\,,
\end{equation}
where $|\cdot|_2$ is the Euclidean norm. Given
$\bm\ve=(\ve_1,\dots,\ve_m)\in(0,+\infty)^m$, the
$\bm\ve$-neighborhood  of the plane $\cL_{\vv a,\vv b}$ is  given by
\begin{equation}\label{e:012}
    \cL^{(\bm\ve)}_{\vv a,\vv b}:=\{\xmn\in\Rmn:|\xmn_j\vv
a+b_j|<\ve_j  \ \ \forall \  j=1, \ldots, m\}\,,
\end{equation}
where $\xmn_j$ is the $j$-th row of $\xmn$. In the case that
$\ve_1=\dots=\ve_m=\ve$, we simply write  $\cL^{(\ve)}_{\vv a,\vv
b}$ for  the symmetric   $\ve$-neighborhood of $\cL_{\vv a,\vv b}$.

\begin{definition}\label{def3}\rm
A non-atomic, finite, doubling Borel measure $\mu$ on $\Rmn$ is
\emph{strongly contracting} if there exist positive constants $C$,
$\alpha$ and $r_0$ such that for any plane $\cL_{\vv a,\vv b}$,
any $\bm\ve=(\ve_1,\dots,\ve_m)\in(0,+\infty)^m$ with $\min \{
\ve_j : 1\le j\le m \} < r_0$ and any $\delta\in(0,1)$ the
following property is satisfied: for all
$\xmn\in\cL^{(\delta\bm\ve)}_{\vv a,\vv b}\cap\Supp$ there is an
open  ball $B$ centred at $\xmn$ such that
\begin{equation}\label{e:013}
B\cap\Supp\subset\cL^{(\bm\ve)}_{\vv a,\vv b} \end{equation} and
\begin{equation}\label{e:014}
\mu(5B\cap\cL^{(\delta\bm\ve)}_{\vv a,\vv b})\le
C\delta^\alpha\mu(5B)\,.
\end{equation}
\end{definition}

\noindent The measure $\mu$ is said to be   \emph{contracting} if
the property holds with $\ve_1=\dots=\ve_m = \ve$.

\bigskip  \noindent\textit{Remark. } The property given by (\ref{e:013}) and (\ref{e:014}) indicates
the rate at which the $\mu$-measure of the $\bm\ve$-neighborhood
of $\cL_{\vv a,\vv b}$ decreases when contracted by the
multiplicative factor $\delta$. Also, note that
\\[0.5ex]
\centerline{$\mu$ is strongly contracting  $ \quad \Longrightarrow
\quad $ $\mu$ is contracting.}

\vspace{2ex}

 The definition of (strongly) contracting is in essence a
global statement -- the `property' is required to hold
for all $\xmn$ in the support $\Supp$.  However, with the view of
establishing `extremality' results such as Theorem \ref{t1}, sets of
$\mu$-measure zero are irrelevant and the notion of (strongly)
contracting almost everywhere suffices. Formally, we say that $\mu$
is \emph{(strongly) contracting almost everywhere} if for
$\mu$-almost every point $\xmn_0 \in \Rmn$ there is a neighborhood
$U$ of $\xmn_0$ such that the restriction $\mu|_{\text{\tiny {\em U
}}}$ of $\mu$ to $U$ is (strongly) contracting.

\subsection{Friendly measures}\label{friendly}

The notion of friendly measures introduced in
\cite{Kleinbock-Lindenstrauss-Weiss-04:MR2134453} identifies purely
geometric conditions on measures on $\R^n$ that are sufficient to
guarantee strong extremality. The class of friendly measures is
defined via the behavior of $\mu$ near hyperplanes $\cL$ in $\R^n$.

\medskip

 Let $\mu$ be a Borel measure on
$\R^{n}$ and as usual let $\Supp $ denote the support of $\mu$. We
say that $\mu$ is \emph{non-planar}\/ if $\mu(\cL)=0$ for any
hyperplane $\cL$. Furthermore, given $\cL$ and a ball $B$ with
$\mu(B)>0$, let $\|d_\cL\|_{\mu,B}$ be the supremum of $\dist(\vv
x,\cL)$ over $\vv x\in\Supp\cap B$. Here $\dist(\vv x,\cL)$ is the
Euclidean distance of $\vv x$ from $\cL$.  Next, let $U$ be an
open subset of $\R^{n}$. Given positive numbers $C$ and $\alpha$,
the measure $\mu$ is called \emph{$(C,\alpha)$-decaying on $U$}\/
if for any non-empty open ball $B\subset U$ centred in $\Supp$,
any affine hyperplane $\cL$ of $\R^{n}$ and any $\ve>0$ one has
that
\begin{equation}\label{e:015}
    \mu(B\cap\cL^{(\ve)})\le
    C\left(\frac{\ve}{\|d_\cL\|_{\mu,B}}\right)^\alpha\,\mu(B)\,.
\end{equation}

\begin{definition}\rm
A non-atomic, Borel measure $\mu$ on $\R^n$ is called
\emph{friendly}\/ if for $\mu$-almost every point $\vv x_0\in\R^n$
there is a neighborhood $U$ of $\vv x_0$ such that  the restriction
$\mu|_{\text{\tiny {\em U }}}$ of $\mu$ to $U$ is finite, doubling,
non-planar and $(C,\alpha)$-decaying for some positive $C$ and
$\alpha$.
\end{definition}

In the next two sections we shall establish that friendly measures
on $\R^n$ identified either with $\R^{1\times n}$ or $\R^{n\times1}$
are strongly contracting.

\subsubsection{Friendly measures on $\R^{m\times1}$}\label{iwdsa}

\begin{proposition}\label{t6}
Any friendly measure $\mu$ on $\R^{m\times 1}$ \!is strongly
contracting almost everywhere.
\end{proposition}

\medskip

\noindent\emph{Proof.} Let $\mu$ be a friendly measure on $\R^{m}$
identified with $\R^{m\times 1}$. Then for $\mu$-almost every point
$\xmn_0 \in \R^{m\times 1} $ there is a neighborhood $U$ of $\xmn_0$
such that $\mu|_U$ is $(C,\alpha)$-decaying on $U$ for some fixed
$C,\alpha>0$. Without loss of generality we can assume that
$\mu=\mu|_U$. The fact that $\mu$ is non-planar, implies that there
are $n+1$ linearly independent points $\xmn_0,\dots,\xmn_n$ in the
support $\Supp$ of $\mu$. This ensures that there exists  a real
number $r_0>0$ such that no single rectangle with shortest  side of
length $\le2r_0$ contains $\Supp$. Also, note that since $\R^m$  is
identified with $\R^{m\times1}$,  the set $\cL_{a,\vv b}$ appearing
in the definition of strongly contracting is simply a point. Thus,
for any point $\cL_{a,\vv b}$ the $\bm\ve$-neighborhood given by
$(\ref{e:012})_{n=1}$ is a rectangle with sides of length $2\ve_i$
($i=1,\dots,m$). In particular, if $\min_{1\le j\le m}\ve_j< r_0$
then
\begin{equation}\label{e:016}
\Supp\not\subset\cL^{(\bm\ve)}_{a,\vv b}\,.
\end{equation}
Fix $\delta \in (0,1)$ and take any point
$\xmn\in\Supp\cap\cL^{(\delta\bm\ve)}_{a,\vv b}$. Without loss of
generality, we  assume that this intersection is non-empty.  The
goal is to construct a ball $B$ centred at $\xmn$ satisfying
(\ref{e:013}) and (\ref{e:014}) in the definition of strongly
contracting. To start with, let $B'$ be an arbitrary ball centred at
$\xmn$ such that
\begin{equation}\label{e:017}
B'\subset\cL^{(\bm\ve)}_{a,\vv b}\,.
\end{equation}
This is possible as $\xmn\in\cL^{(\bm\ve)}_{a,\vv b}$ and by
definition $\cL^{(\bm\ve)}_{a,\vv b}$ is an open set. By
(\ref{e:016}) and (\ref{e:017}), there is a real number $\tau\ge1$
such that
\begin{equation}\label{e:018}
5\tau B'\cap\Supp\not\subset\cL^{(\bm\ve)}_{a,\vv b} \qqand \tau
B'\cap\Supp\subset\cL^{(\bm\ve)}_{a,\vv b}\,.
\end{equation}
By (\ref{e:018}), there exists a point $ \xmn'\in \Big(5\tau
B'\cap\Supp\Big)\setminus\cL^{(\bm\ve)}_{a,\vv b}\,. $ By the choice
of $\xmn'$, there exists a $j\in \{1,\dots,m\}$ such that
\begin{equation}\label{e:019}
\big|X_{j}' +b_j\big|\ge \ve_j\,.
\end{equation}

\noindent Recall that $X_{j}'$ and $ b_j $ are the $j$-th
coordinates of $\xmn$ and $\vv b$ respectively. With reference to
\S\,\ref{friendly}, let $\cL$ be given by $X_j+b_j=0$ and $B=5\tau
B'$. It follows from (\ref{e:019}) that $ \|d_{\cL}\|_{\mu,B}\ge
\ve_j\,. $ Since $\mu$ is $(C,\alpha)$-decaying, (\ref{e:015}) with
$\ve:=\delta\ve_j$ implies that
$$
\mu(5\tau B'\cap\cL^{(\delta\bm\ve)}_{a,\vv b}) \ \le \ \mu(5\tau
B'\cap \cL^{(\delta\ve_j)}) \ < \
C\left(\frac{\delta\ve_j}{\ve_j}\right)^\alpha\mu(5\tau B') \ = \
C\,\delta^\alpha\mu(5\tau B')\,.
$$
The upshot of this is that the ball $\tau B'$ satisfies conditions
(\ref{e:013}) and (\ref{e:014}). The other conditions of strongly
contracting are trivially met and the proof is complete.

\hfill $\boxtimes$

\subsubsection{Friendly measures on $\R^{1\times n}$}

\begin{proposition}\label{t7}
Any friendly measure $\mu$ on $\R^{1\times n}$ is strongly
contracting almost everywhere.
\end{proposition}

\medskip

\noindent\emph{Proof.}  Let $\mu$ be a friendly measure on $\R^{n}$
identified with $\R^{1\times n}$. Then for $\mu$-almost every point
$\xmn_0 \in \R^{1\times n}$ there is a neighborhood $U$ of $\xmn_0$
such that $\mu|_U$ is $(C,\alpha)$-decaying on $U$ for some fixed
constants $C,\alpha>0$. Without loss of generality we can assume
that $\mu=\mu|_U$. The fact that $\mu$ is non-planar, implies that
there are $n+1$ linearly independent points $\xmn_0,\dots,\xmn_n$ in
the support $\Supp$ of $\mu$. This ensures that there exists  a real
number  $r_0>0$ such that for any $\ve \in (0,r_0)$  and any
hyperplane $\cL_{\vv a,b}$ the $\ve$-neighborhood given by
$(\ref{e:012})_{m=1}$  cannot contain all the points
$\xmn_0,\dots,\xmn_n$. It follows that for any hyperplane and  $0<
\ve<r_0$, we have that
\begin{equation}\label{e:020}
\Supp\not\subset\cL_{\vv a,b}^{(\ve)}   \ .
\end{equation}
Fix $\delta \in (0,1)$ and take any point
$\xmn\in\Supp\cap\cL^{(\delta\bm\ve)}_{a,\vv b}$ -- we may as well
assume that this intersection is non-empty.  Now, let $B'$ be an
arbitrary ball centred at $\xmn$ such that
\begin{equation}\label{e:021}
B'\subset\cL_{\vv a,b}^{(\ve)}\,.
\end{equation}
This is possible as $\cL_{\vv a,b}^{(\ve)}$ is open. By
(\ref{e:020}) and (\ref{e:021}), there is a real number $\tau\ge 1$
such that
\begin{equation}\label{e:022}
5\tau B'\cap\Supp\not\subset\cL_{\vv a,b}^{(\ve)} \qqand \tau
B'\cap\Supp\subset\cL_{\vv a,b}^{(\ve)}\,.
\end{equation}
By (\ref{e:022}), there exists  a point $\xmn'\in5\tau B'\cap\Supp$
which is not contained in $\cL_{\vv a,b}^{(\ve)}$. With reference to
\S\,\ref{friendly}, let $\cL=\cL_{\vv a,b}$ and $B=5\tau B_0$. It
follows that  $ \|d_{\cL}\|_{\mu,B}\ge \ve\, $ which together with
(\ref{e:015}) implies that
$$   \mu\big(5\tau B'\cap\cL_{\vv a,b}^{(\delta\ve)}\big)
    \ \le \   C\ \delta^\alpha \ \mu\big(5\tau B'\big)\,.
$$
Thus, the ball $\tau B'$ satisfies conditions (\ref{e:013}) and
(\ref{e:014}) in the definition of strongly contracting. The other
conditions of strongly contracting are trivially met and the proof
is complete.

\hfill $\boxtimes$

\subsection{Differentiable manifolds: Proof of Theorem~\ref{t4}\label{svdiffman}}

We begin by establishing  Theorem~\ref{t4}B - part (B) of Theorem
\ref{t4}. Clearly, we only have to prove the necessity part as the
right hand side of the statement contains the left hand side. Thus,
we are given that the differentiable submanifold $\cM$ of $\R^n$ is
strongly extremal. Let $m$ denote the Riemannian measure on $\cM$.
The aim is to show that $m$ is inhomogeneously strongly extremal on
$\R^{n\times1}$ -- see Definition~\ref{def2} in \S\ref{inhomtheory}.
This is a simple consequence of Theorem \ref{t1}B once we have
established that $m$ as a measure on $\R^{n\times1}$ is  strongly
contracting almost everywhere.

Take any point $\vv y_0\in\cM$ such that the tangent plane to $\cM$
at $\vv y_0$ is not orthogonal to any of the coordinate axes. Since
the latter property holds almost everywhere on $\cM$, it suffices to
prove Theorem \ref{t4}B for a neighborhood $\cP$ of $\vv y_0$.
Without loss of generality, we can assume that there is a $C^{(1)}$
local parameterisation of $\cP$ given by $\vv f:U\to\cM$. Here $U$
is a ball in $\R^d$ centred at $\vv x_0\in U$ such that $\vv f(\vv
x_0)=\vv y_0$ and $d=\dim\cM$. By the condition on $\vv y_0$ imposed
above, there is a direction $\vv v\in\R^d$ such that the tangent
direction $\frac{\partial \vv f(\vv x_0)}{\partial \vv v}$ is not
orthogonal to any of the coordinate axes. This means that there
exists some $\kappa>0$, such  that
$$ 2\kappa^{-1}<\left|\frac{\partial f_i(\vv
x_0)}{\partial \vv v}\right|<\kappa/2 \qquad\text{for all $1\le i\le
n$.}
$$
 Since $\vv f$ is $C^{(1)}$, there exists  a sufficiently
small ball $B_0\subset U$ centered at $\vv x_0$ such that
\begin{equation}\label{e:023}
 \kappa^{-1}<\left|\frac{\partial f_i(\vv x)}{\partial \vv
 v}\right|<\kappa\qquad\text{for all $1\le i\le n$ and all $\vv x\in B_0$.}
\end{equation}

\noindent  Without loss of generality, take $f(B_0)$ to be the
neighborhood $\cP \subset \cM$ of $\vv y_0$ mentioned above. We now
slice $B_0$ with respect to the direction $\vv v$ so as to reduce
the problem at hand to one concerning differentiable curves. Since
$\cM$ is strongly extremal and using the fact that sets of full
measure are invariant under diffeomorphisms, the set
$$\cE:=\{\vv x\in B_0:\ws(\vv f(\vv x))=1\}$$
has full Lebesgue measure in $B_0$. Now for any $\vv x'\in\R^{d}$
orthogonal to $\vv v$, consider the line $L_{\vv x'}$ in $\R^d$
given by
$$
L_{\vv x'}\,:=\{\vv x=x\vv v+\vv x'\in\R^d: x\in\R\}\,.
$$
Also,  let
$$\cE_{\vv x'}:=\cE\cap L_{\vv x'} \qquad \text{ and  } \qquad B_{\vv
x'}:=B_0\cap L_{\vv x'}  \ .  $$ Clearly, $B_{\vv x'}$ is either an
interval or is empty and  $\cE_{\vv x'}\subset B_{\vv x'}$. For
obvious reasons, we only consider the situation when $B_{\vv
x'}\not=\emptyset$. Since $\cE$ has full measure in $B_0$, it
follows from Fubini's theorem that for almost every $\vv x'$ the
slice $\cE_{\vv x'}$ has full measure in $B_{\vv x'}$.  Now let $\vv
f_{\vv x'}$ denote the map   $\vv f$ restricted to $B_{\vv x'}$.
Clearly,  $\vv f_{\vv x'}$ is  a diffeomorphism from $B_{\vv x'}$
onto the curve
$$\cM_{\vv x'}:=\vv f(B_{\vv x'}) \, . $$ Since
$\cE_{\vv x'}$ has full measure in $B_{\vv x'}$ and $\vv f_{\vv x'}$
is a diffeomorphism, $\cM_{\vv x'}$ is strongly extremal for almost
all $\vv x'$ orthogonal to the direction $\vv v$.

\noindent Now we fix any $\vv x'\in\R^{d}$ orthogonal to $\vv v$
such that the curve $\cC:=\cM_{\vv x'}$ is non-empty and strongly
extremal. Define the map $\vv g=(g_1,\dots,g_n):I\to\R^n$ from the
interval
$$ I:=\{x\in\R:x\vv v+\vv x'\in B_0\}$$ such that $\vv
g(x)=\vv f(x\vv v+\vv x')$. By (\ref{e:023}), we have that
\begin{equation}\label{e:024}
    \kappa^{-1}\le |g_i'(x)|\le \kappa  \qquad \text{for all $1\le i\le n$ and all  $x\in I$.}
\end{equation}

\noindent Let $\mu$ denote the induced Lebesgue measure on $\cC$ and
identify $\R^n$ with $\R^{n\times1}$. The key part of the proof is
to that show that $\mu$ as a measure on $\R^{n\times1}$ is strongly
contracting. This involves verifying (\ref{e:013}) and
(\ref{e:014}). Regarding (\ref{e:014}) we can assume that
$\delta<1/2$ as otherwise (\ref{e:014}) is trivially satisfied with
$C>2$. Also, given that $\R^n$  is being identified with
$\R^{n\times1}$,  the set $\cL_{a,\vv b}$ appearing  in the
definition of strongly contracting is simply a  point. Now choose a
real number $r_0>0$ so that for any point  $\cL_{a,\vv b}$ and
$\bm\ve=(\ve_1,\dots,\ve_n)\in(0,+\infty)^n$ with $\min_{1\le j\le
n}\ve_j< r_0$, we have that
$$
\cC\not\subset\cL_{a, \vv b}^{(\bm\ve)} \ .
$$
The latter is readily deduced from (\ref{e:024}). In what follows
we fix a point $\cL_{a,\vv b}$, a vector $\bm\ve$ with $0<
\min_{1\le j\le n}\ve_j< r_0$ and a $\delta\in(0,1/2)$. Without
loss of generality, we assume that $\ve_1\le\dots\le\ve_n$ and
that $\mu (\cL_{a,\vv b}^{(\delta\bm\ve)}\cap\cC ) \not= 0$. We
now verify that
\begin{equation}\label{e:025}
 \mu(\cL_{a,\vv b}^{(\delta\bm\ve)}\cap\cC) \, \le  \,2\sqrt
n\kappa^2\delta\ve_1\,.
\end{equation}
 Let $\xmn$ and $\xmn'$ be any  two points in
$\cL_{a,\vv b}^{(\delta\bm\ve)}\cap\cC$. Thus, $\xmn=\vv g(x)$ and
$\xmn'=\vv g(x')$  for some $x,x'\in I$. It follows that
\begin{equation}\label{e:026}
|g_1(x)-b_1/a|<\delta\ve_1\qquad\text{and}\qquad|g_1(x')-b_1/a|<\delta\ve_1\,,
\end{equation}
where $b_1$ is the first coordinate of $\vv b$ associated with the
point $\cL_{a,\vv b}$ and $g_1$ is the first coordinate function of
$\vv g$. By the Mean Value theorem, there exists some
$\theta\in[0,1]$ such that
\begin{eqnarray}\label{e:027}
|(g_1(x)-b_1/a)-(g_1(x')-b_1/a)|& = & |g_1(x)-g_1(x')|  \nonumber \\
&=&|x-x'|\,\big|g_1'\big(\theta x+(1-\theta')x'\big)\big|   \ .
\end{eqnarray}

By (\ref{e:024}), (\ref{e:026}) and (\ref{e:027}), it follows that
\begin{equation}\label{e:028}
    |x-x'| \, \le \, 2\kappa\delta\ve_1 \ .
\end{equation}
In view of  (\ref{e:024}), $g_i$ is monotonic for every $i$ and
therefore the set $\vv g^{-1}(\cL_{a,\vv
b}^{(\delta\bm\ve)}\cap\cC)\subset I$ is an interval. Let $x_0 $ and
$x_1$ be the endpoints of this interval with $x_0<x_1$. Clearly,
(\ref{e:028}) is valid  with $x=x_0$ and $x'=x_1$. Therefore,
$$
\mu(\cL_{a,\vv b}^{(\delta\bm\ve)}\cap\cC)=\int_{x_0}^{x_1}|\vv
g'(x)|dx\stackrel{(\ref{e:024})}\le \sqrt
n\kappa|x_0-x_1|\stackrel{(\ref{e:028})}{\le} \sqrt n\kappa
2\kappa\delta\ve_1\,.
$$
This is precisely (\ref{e:025}). Now, let $B$ be a ball centred at
$\xmn$ of radius $\ve_1/2$. By the choice of $\xmn$ and the fact
that $\delta<1/2$, we have that
\begin{equation}\label{e:029}
B\subset\cL_{a,\vv b}^{(\bm\ve)}.
\end{equation}
By the  Mean Value Theorem and (\ref{e:024}),  for any $x\in
I':=\{x'\in I:|x-x'|<\ve_1/(2\sqrt n\kappa)\}$ we have that
$$
|\vv g(x')-\vv g(x)|=|\vv g'(\theta' x'+(1-\theta')x)|\
|x-x'|<\ve_1/2\,.
$$
It follows that $\vv g(x')\in B$ for any $x'\in I'$ and that
$|I'|\ge \ve_1/(2\sqrt n\kappa)$. Hence,
\begin{equation}\label{e:030}
\mu(5B\cap\cC)\ge \mu(B\cap\cC)\ge\int_{I'}|\vv
g'(x)|dx\stackrel{(\ref{e:024})}{\ge}\kappa\sqrt
n/2\cdot|I'|\ge\ve_1/4\,.
\end{equation}
On combining inequalities (\ref{e:025}) and (\ref{e:030}), we obtain
that
\begin{equation}\label{e:031}
\mu(5B\cap\cL_{a,\vv b}^{(\delta\bm\ve)} \cap \cC)\le \mu(\cL_{a,\vv
b}^{(\delta\bm\ve)}\cap\cC)\le 2\sqrt n\kappa^2\delta\ve_1\le 8\sqrt
n\kappa^2\delta \, \mu(5B \cap \cC)\,.
\end{equation}
Clearly, (\ref{e:029}) verifies (\ref{e:013}) and (\ref{e:031})
verifies (\ref{e:014}). Thus, the measure $\mu$ on $\R^{n \times 1}
$  is strongly contracting.  By Theorem~\ref{t1}B, it follows that
$\mu$ is inhomogeneously strongly extremal on $\R^{n\times1}$. By
definition, $\mu$ or equivalent $\cC$ is simultaneously
inhomogeneously strongly extremal. This establishes Theorem
\ref{t4}B in the case that $\cM$ is a differentiable curve.
\noindent To deal with manifolds in general, we appeal to Fubini's
theorem.  For any $\bm\theta\in\R^n$, consider the sets
$$
 \cE^{\bm\theta}:=\{\vv x\in B_0:\ws(\vv f(\vv x),\bm\theta)=1\}
\qquad \text{ and  } \qquad \cE_{\vv
x'}^{\bm\theta}:=\cE^{\bm\theta}\cap L_{\vv x'}  \ . $$
 Clearly,
$\cE^{\bm\theta}_{\vv x'}\subset B_{\vv x'}$. For almost every $\vv
x'$ the measure $\mu$ on the corresponding curve $\cM_{\vv x'}$ is
simultaneously inhomogeneously strongly extremal. Thus, for almost
every $\vv x'$ the slice  $\cE^{\bm\theta}_{\vv x'}$ has  full
Lebesgue measure in $B_{\vv x'}$. Hence, by Fubini's theorem we have
that $\cE^{\bm\theta}$ has full Lebesgue measure in $B_0$.
Consequently, $\vv f(\cE^{\bm\theta})$ has full Riemannian measure
in $\cP := \vv f(B_0)$. This completes the proof of
Theorem~\ref{t4}B.

\hfill $\boxtimes$

\medskip

The proof of Theorem~\ref{t4}A follows the same line of argument as
above. However, we only require that the inequality  in
(\ref{e:023}) holds for at least one value of $i$ rather than for
all $i$. This is the case for any differentiable manifold
irrespective of the direction $\vv v$. Hence there is no extra
hypothesis  on $\cM$ in Theorem~\ref{t4}A. The details are left to
the reader. As mentioned at the end of \S\ref{svresults}, for a
self-contained and independent proof of Theorem \ref{t4}A see
\cite{Beresnevich-Velani-Moscow}.

\section{Lower bounds for Diophantine exponents\label{TI} }

Given a measure $\mu$ on $\Rmn$, suppose we are interested in
establishing that $\mu$ is inhomogeneously strongly extremal.
Clearly, this would follow on showing that for all $\bth\in\R^m$
$$
\wsxth \le 1\quad\text{for $\mu$-almost all  $\xmn \in \Rmn $}
$$
and
\begin{equation}\label{e:032}
\wsxth \ge 1\quad\text{for $\mu$-almost all   $\xmn \in \Rmn$}.
\end{equation}

\noindent Establishing inhomogeneous extremality corresponds to
similar statements with $\wsxth$ replaced by $\wxth$. Note that as a
consequence of (\ref{e:003}) and (\ref{e:004}),  in the homogeneous
case ($\bth=\vv0$) the set of $\xmn$ satisfying $\wxth \ge 1$ or
$\wsxth \ge 1$ is  the whole space. Thus, (\ref{e:032}) is
automatically satisfied within the homogeneous setting. A priori,
this is not the case within the inhomogeneous setting. The goal of
this section is to establish (\ref{e:032}) and the analogous  $\wxth
\ge 1 $ statement for extremal measures within the inhomogeneous
setting.

\begin{proposition}\label{lem1}
Let $\mu$ be an  extremal measure  on $\Rmn$. Then for all
$\bth\in\R^{m}$,
$$
\wsxth \ge \wxth \ge 1\quad\text{for $\mu$-almost all  $\xmn \in
\Rmn $} \, .
$$
\end{proposition}

\noindent If $\mu$ is extremal,  then $\wx=1$ for $\mu$-almost all
$\vv X\in\Rmn$ and Proposition~\ref{lem1} readily follows from
(\ref{e:004}) and the following statement.

\begin{lemma}\label{lem2}
Let $\xmn\in\Rmn$ such that $\wx=1$. Then for all $\bth\in\R^m$,
\begin{equation}\label{e:033}
\wxth\ge1\,.
\end{equation}
\end{lemma}

The proof of the lemma utilises  basic `transference' inequalities
relating various forms of Diophantine exponents. These we briefly
describe. A form of Khintchine's transference principle due to Dyson
\cite[Theorem~5C]{Schmidt-1980} relates the homogeneous exponents of
$\xmn$ and its transpose $\txmn$. It states that
\begin{equation}\label{e:034}
 \wx=1 \iff w(\txmn)=1   \qquad\text{for all }\xmn\in\Rmn \,.
\end{equation}
In the spirit of  Cassels \cite[Chapter~5]{Cassels-1957}, Bugeaud \&
Laurent \cite{Bugeaud-Laurent2005} have recently discovered
transference inequalities that relate the Diophantine exponents
$\wxth$ with their uniform counterparts $\hwxth$. The latter are
defined as followed. Given $\xmn\in\Rmn$ and $\bth \in\R^{m}$, let
$\hwxth$ be the supremum of $w\ge0$ such that for all sufficiently
large $Q$ there is a $\vv q\in\Z^n\bnz$ satisfying (\ref{e:001}). As
with the standard non-uniform exponents, a trivial consequence of
Dirichlet's theorem is that
\begin{equation}\label{e:035}
 \hwx\ge 1\qquad\text{for all }\xmn\in\Rmn\,.
\end{equation}
Also, the following inequalities are easily verified.
\begin{equation}\label{e:036}
 \wxth \ \ge \ \hwxth\ \ge \ 0\,.
\end{equation}

\bigskip

\noindent\textbf{Theorem BL (Bugeaud \& Laurent) } {\it Let
$\xmn\in\Rmn$. Then for all $\bth \in\R^m$,
\begin{equation}\label{e:037}
 \wxth\ge\frac1{\wh w(\txmn)}
 \qqand
 \hwxth\ge\frac1{w(\txmn)}
\end{equation}
with equalities in $(\ref{e:037})$ for almost all $\bth\in\R^m$. }

\bigskip

We are now fully armed to proceed with the  proof of above lemma.

\medskip

\noindent\emph{Proof of Lemma~\ref{lem2}.}  We are given that
$\wx=1$. Hence, by (\ref{e:034}) it follows that $w(\txmn)=1$.
This together with  (\ref{e:035}) and
(\ref{e:036})${}_{\bth=\vv0}$ applied to $\txmn$ implies that $\wh
w(\txmn)=1$. In turn, this combined with (\ref{e:037}) implies
that $\wxth\ge1$ and thereby completes the proof.

\vspace{-2ex}

\hfill $\boxtimes$

\bigskip

\noindent {\em Remark 1. } It is worth pointing out that
Lemma~\ref{lem2}, which allows us to deduce Proposition \ref{lem1}
and thereby reduce the proof of Theorem \ref{t1} to establishing
upper bounds for the associated Diophantine exponents, can in fact
be proved without appealing to Theorem BL. Indeed, a proof can be
given which only makes use of classical transference inequalities;
namely Theorem VI of Chapter 5 in \cite{Cassels-1957}. Thus, the
proof of Proposition \ref{lem1} and therefore  Theorem \ref{t1} is
not actually dependent on the recent developments regarding
transference inequalities.

\medskip

\noindent {\em Remark 2. } Theorem BL actually gives us information
beyond Lemma~\ref{lem2}. It enables us to deduce that inequality
(\ref{e:033}) is in fact an equality for almost all $\bth \in\R^m$.
Thus, the real significance of Theorem~\ref{t1}A is  in establishing
a global result which holds for \emph{all}\/ $\bth\in\R^m$.

\subsection{Proof of Proposition~\ref{prop1}} \label{sec:1B}

Let $\mu$ be a measure on $\Rmn$. Given that $\mu$ is
inhomogeneously strongly extremal we wish to conclude that $\mu$ is
inhomogeneously  extremal. This as we shall now see is a simple
consequence of (\ref{e:004}) and Lemma~\ref{lem2}.

 We are given that for any $\bth\in\R^m$, $\wsxth=1$  for
$\mu$-almost all $\xmn\in\Rmn$. By (\ref{e:004}), it follows that
for any $\bth\in\R^m$, $\wxth\le1$ for $\mu$-almost all
$\xmn\in\Rmn$. Thus, we only need to show that for any
$\bth\in\R^m$, $\wxth\ge1$ for $\mu$-almost
all $\xmn\in\Rmn$. Since $\mu$ is inhomogeneously strongly extremal
we trivially have that $\mu$ is strongly extremal and therefore
extremal. In other words,  $\wx=1$ for almost all $\xmn\in\Rmn$.
This together with Lemma~\ref{lem2} yields the desired statement.

\vspace{-2ex}

\hfill $\boxtimes$

\section{A reformulation of Theorem \ref{t1}}\label{sec:1}

The goal of this section is to reformulate Theorem \ref{t1} so that
the new statement can be deduced via the general framework developed
in \S\ref{HIT}. On the other hand, the reformulation is natural even
for a direct proof of Theorem \ref{t1} and thereby motivates the
general framework.

Theorem~\ref{t1} consists of two parts which we refer to as Theorem
\ref{t1}A and Theorem \ref{t1}B. We will concentrate on establishing
Theorem \ref{t1}B. The proof of Theorem~\ref{t1}A is similar in
spirit and we shall indicate the necessary modifications that need
to be made.

\noindent With the intention of proving Theorem \ref{t1}B, let $\mu$
be a strongly extremal measure  on $\Rmn$ and define
$$\cA_{m,n}^{\bth}:=\{\xmn\in\Rmn:\wsxth>1\} \,   . $$

\noindent In view of Proposition~\ref{lem1}, Theorem \ref{t1}B is
reduced to showing  that
\begin{equation}\label{e:038}
\mu(\cA_{m,n}^{\bth})=0  \qquad\text{for all }\bth\in\R^m \,.
\end{equation}

\noindent The key towards establishing (\ref{e:038}) is the
following reformulation. Let $\TTT$ denote a countable subset of
$\R^{m+n}$ such that for every $\ttt=(t_1,\dots,t_{m+n})\in\TTT$

\begin{equation}\label{e:039}
    \sum_{j=1}^mt_j=\sum_{i=1}^{n}t_{m+i}\,.
\end{equation}

\noindent For $\ttt\in\TTT$, consider the diagonal unimodular
transformation $g_{\vv t} $ of $\R^{m+n}$ given by
\begin{equation}\label{e:040}
  g_{\vv t} \ := \ \operatorname{diag}\{2^{t_1},\dots,2^{t_m},2^{-t_{m+1}},\dots,2^{-t_{m+n}}\}\,.
\end{equation}

\noindent For $\xmn \in \Rmn$, define the matrix
$$
 M_{\xmn}:=\left(\begin{array}{cc}
  I_m&\xmn\\[2ex]
  0&I_n
\end{array}
 \right)\,,
$$
where $I_n$ and $I_m$ are respectively the  $n\times n$ and $m\times
m$ identity matrices. The matrix $M_{\xmn}$ is a linear
transformation of $\R^{m+n}$. Given $\bth \in \R^m$, let
$$
M_{\xmn}^{\bth}\ :\ \vv a\mapsto M_{\xmn}^{\bth}\vv a:=M_{\xmn}\vv
a+\bm\Theta\,,
$$

\noindent where $\bm\Theta \, := \,
{}^t(\theta_1,\dots,\theta_m,0,\dots,0)\in\R^{m+n}$. Thus,
$M_{\xmn}^{\bth}$ is  an affine transformation of $\R^{m+n}$.

\noindent Let
\begin{equation}\label{e:041}
\AAA=\Z^m\times(\Z^n\bnz)\,.
\end{equation}
Then, for  $\ve>0$, $\ttt\in\TTT$ and $\alpha\in\cA$ define the sets
\begin{equation}\label{e:042}
\De^{\bth}_\ttt(\alpha,\ve):=\{\xmn\in\Rmn:| g_{\vv t}M_{\vv
X}^{\bth}\alpha|\,< \ve\}
\end{equation}
and
$$
\De^{\bth}_\ttt(\ve):=\bigcup_{\alpha\in
\AAA}\De^{\bth}_\ttt(\alpha,\ve)= \{\xmn\in\Rmn:\inf_{\alpha\in\cA}|
g_{\vv t}M_{\xmn}^{\bth}\alpha|\,< \ve\}\,.
$$

\noindent For $\eta>0$,  define the function
\begin{equation}\label{e:043}
\psi^\eta \, : \, \TTT \mapsto \R_+   \  :  \
\ttt\mapsto\psi^\eta_\ttt :=2^{-\eta\sigma(\ttt)} \
\end{equation}
where $\sigma(\ttt):=t_1+\dots+t_{m+n}$, and consider the $\limsup$
set given by
\begin{equation}\label{e:044}
 \La^{\bth}_{\TTT}(\psi^\eta\,) \, := \, \limsup_{\ttt \in \TTT }\De^{\bth}_\ttt(\p^\eta_\ttt)\, .
\end{equation}
In the case $\bth=\vv0$, we write $\La_{\TTT}(\p^\eta)$ for
$\La^{\bth}_{\TTT}(\p^\eta)$. The following result provides a
reformulation of the set $ \cA_{m,n}^{\bth} $ in terms of the
$\limsup$ sets given by  (\ref{e:044}).

\begin{proposition}\label{prop5}
There exists a countable subset $\TTT$ of $\R^{m+n}$ satisfying
$(\ref{e:039})$ such that
\begin{equation}\label{e:045}
 \sum_{\ttt\in\TTT} 2^{- \eta\, \s t} < \ \infty \   \  \qquad  \forall \ \eta > 0
\end{equation}
and
\begin{equation}\label{e:046} \cA_{m,n}^{\bth}= \bigcup_{\eta >
0}\La^{\bth}_{\TTT}(\psi^\eta\,) \  \qquad \forall \ \bth\in\R^m \,
. \end{equation}
\end{proposition}

\bigskip

\noindent  Now,  let $\mu$ be a measure on $\Rmn$ that is strongly
contracting. Then,  as a consequence of (\ref{e:038}) and
Proposition \ref{prop5}, the proof of Theorem~\ref{t1}B is reduced
to showing that
\begin{equation}\label{e:047}
       \mu(\La_{\TTT}(\p^\eta))=0     \  \quad  \forall \ \eta > 0   \qquad\Longrightarrow\qquad
       \mu(\La^{\bth}_{\TTT}(\p^\eta))=0 \  \quad  \forall \ \eta > 0  \ \ \ \ \
\end{equation}
for a suitable choice of $\TTT$ satisfying $(\ref{e:039})$ and
$(\ref{e:045})$. It is worth mentioning, that in the proof of
Proposition \ref{prop5} an  explicit choice of $\TTT$ is given.

\vspace{1ex}

\noindent {\em Remark. \ } In the case of Theorem~\ref{t1}A, the
analogous reformulation and reduction are equally valid. The only
difference is that $\cA_{m,n}^{\bth} $ is defined in terms of $
\wxth $ rather than $\wsxth $.

\subsection{Proof of Proposition~\ref{prop5}}

Given  $\vv s =(s_1,\dots,s_{m}) \in\Z_+^m$ and $\vv l
=(l_1,\dots,l_{n}) \in\Z_+^n$, let
$$ \s s \, := \,
\sum_{j=1}^ms_j   \ \  , \ \ \quad \s l \, := \, \sum_{i=1}^nl_i
\qquad \text{and} \qquad \zeta:=\zeta(\vv s,\vv l)=\dfrac{\s s-\s
l}{m+n} \ ,
$$
where $\Z_+$ is the set of non-negative integers. Furthermore,
define the $(m+n)$-tuple $\vv t=(t_1,\dots,t_{m+n})$ by setting
\begin{equation}\label{e:048}
    \vv t:=\Big(s_1-\zeta,\dots,\ s_m-\zeta,\ l_1+\zeta,
    \dots,l_n+\zeta\Big)\,  \
\end{equation}
\noindent  and let
\begin{equation}\label{e:049}
    \TTT:=\{\ttt\in\R^{m+n}\text{ defined by (\ref{e:048})}:\vv s\in\Z_+^m,\,\vv
    l\in\Z_+^n \text{ with } \s s \ge \s t \}\, .
\end{equation}

\medskip

\noindent The goal is to show that this choice of $\TTT$ is
suitable within the context of  Proposition \ref{prop5}. Equality
(\ref{e:039}) readily follows from (\ref{e:048}). Next, it is
easily verified that for any $ \ttt \in \TTT$
\begin{equation}\label{e:050}
    \tfrac12\ \s t=\s s-m \zeta=\s l+n \zeta\,,
\end{equation}
where $
 \s t \, := \,
\sum_{k=1}^{m+n}t_k \, $. This together with the fact that $
\zeta$ is non-negative, yields that
\begin{equation}\label{e:051}
\s l \le \tfrac12\,\s t \le \s s \, .
\end{equation}
Furthermore, on summing the two different expressions for
$\frac12\s t$ arising in (\ref{e:050}) and using the fact that $\s
l\ge0$, we obtain that
\begin{equation}\label{e:052}
\begin{array}[b]{rcl}
\s t  = \s s+\s l-\frac{m-n}{m+n}\big(\s s-\s l\big)
   & \ge &  \s s+\s l-\frac{|m-n|}{m+n}\big(\s s+\s l\big) \\[2ex]
    & \ge &  \frac{1}{m+n}\big(\s s+\s l\big) \; .
\end{array}
\end{equation}
The latter inequality establishes (\ref{e:045}). In turn, it follows
that for any $v\in\R_+$
\begin{equation}\label{e:053}
\#\{\, \ttt\in\TTT\, :\, \s t\, < v\, \}\, < \, \infty  \ .
\end{equation}

\medskip

\noindent Now to establish the set equality (\ref{e:046}), fix $
\bth\in\R^m $.   It is easily verified that
$\xmn\in\cA^{\bth}_{m,n}$ if and only if there exists  an $\ve>0$,
such that for arbitrarily large $Q>1$ there is an $\alpha=(\vv
p,\vv q)\in\cA:=\Z^m\times(\Z^n\bnz)$ satisfying $|\xmn\vv q+\vv
p+\bth|\le1/2$ such that
\begin{equation}\label{e:054}
 \PI\big(\xmn\vv q+\vv p+\bth\big)<Q^{-(1+\ve)}\qqand\PI_+(\vv q)\le
 Q \ .
\end{equation}

\bigskip

\noindent {\bf Step 1. } We  show that
\begin{equation}\label{e:055}
\cA_{m,n}^{\bth}\subseteq \bigcup_{\eta >
0}\La^{\bth}_{\TTT}(\psi^\eta\,) \, .
\end{equation}

\noindent Suppose $\xmn\in\cA^{\bth}_{m,n}$.  It follows that
(\ref{e:054}) is satisfied for infinitely many $Q\in\Z_+$.  For any
such $Q$, there exist unique $\vv s\in\Z_+^m $ and $\vv l\in\Z_+^n\
$ such that
\begin{equation}\label{e:056}
    2^{-s_j}\le \max\Big\{|\xmn_{j} \, \vv q+p_j+\theta_j|\ ,\
    Q^{-(1+\ve)}\Big\}<2^{-s_j+1}\qquad\text{for } \ \ 1\le j\le m
\end{equation}
and
\begin{equation}\label{e:057}
    2^{\, l_i}\le \max\{1,|q_i|\}<2^{\, l_i+1}\qquad\text{for } \ \ 1\le i\le
    n\,.
\end{equation}

\noindent Throughout, $\xmn_{j}:=(x_{j,1},\dots,x_{j,n})$ denotes
the $j$-th row of $\xmn\in\Rmn$. By (\ref{e:056}) and (\ref{e:057}),
we have that
$$2^{\s l}\le\PI_+(\vv q) \qquad \text{ and } \qquad
2^{-\s s}<\max\Big\{\PI \big(\xmn\vv q+\vv p+\bth\big) \ , \
Q^{-(1+\ve)}\Big\} \, . $$ This together with  (\ref{e:054}) implies
that $ 2^{-\s s}<2^{-\s l(1+\ve)}$. Hence,
\begin{equation}\label{e:058}
\s s-\s l>\ve \s l\ge0\,.
\end{equation}
Thus, $ \ttt $ given by (\ref{e:048}) with  $ \vv s =
(s_1,\dots,s_{m})$ and $\vv l = (l_1,\dots,l_{n})$ satisfying
(\ref{e:056}) and (\ref{e:057}) is  in  $\TTT$.

\noindent  If $\s s\le 2\s l$, then
$$
\zeta \, = \, \frac{\s s-\s l}{m+n} \, \stackrel{(\ref{e:058})}{\ge}
\,  \frac{\ve \s l}{m+n} \, \ge \,  \frac{\ve \s s}{2(m+n)}  \,
\stackrel{(\ref{e:051})}{>} \, \frac{\ve \s t}{4(m+n)} \ .
$$
If $\s s>2\s l$,  then
$$
\zeta \, = \, \frac{\s s-\s l}{m+n} \, \ge \,  \frac{\s s}{2(m+n)}
\, \stackrel{(\ref{e:051})}{>}  \, \frac{\s t}{4(m+n)}\ .
$$
On combining the above inequalities,  we deduce  that
\begin{equation}\label{e:059}
    \zeta \, > \, \eta_0\; \s t\qquad\text{with } \quad \eta_0:=\frac{\min\{\ve,1\}}{4(m+n)}\,.
\end{equation}

\bigskip

\noindent With reference to (\ref{e:040}), we have that
$$
  g_{\vv t} \ = \
  2^{-\zeta}\operatorname{diag}\{2^{s_1},\dots,2^{s_m},2^{-l_1},\dots,2^{-l_n}\}\,
$$
and in view of  (\ref{e:056}) and (\ref{e:057}), it follows  that
\begin{equation}\label{e:060}
\inf_{\alpha\in\cA}| g_{\vv t}M_{\xmn}^{\bth}\alpha|\,< 2\cdot
2^{-\zeta}\,.
\end{equation}
For $0<\eta<\eta_0$, (\ref{e:059}) together with (\ref{e:060})
implies that
\begin{equation}\label{e:061}
\inf_{\alpha\in\cA}| g_{\vv t}M_{\xmn}^{\bth}\alpha|\,< 2^{-\eta\,
\s t}
\end{equation}
for all sufficiently large $\s t$. Note that (\ref{e:054}) and
(\ref{e:056}) ensure that $\s s\to\infty$ as $Q\to\infty$.
Therefore, in view of (\ref{e:052}) and the fact that (\ref{e:054})
is satisfied for infinitely many $Q\in\Z_+$, we have that
(\ref{e:061}) is satisfied for infinitely many $ \ttt \in \TTT$. The
upshot is that $\xmn \in \La^{\bth}_{\TTT}(\psi^\eta\,) $ for any
$\eta \in (0, \eta_0) $. This establishes (\ref{e:055}).

\bigskip

\noindent {\bf Step 2. } We show that
\begin{equation}\label{e:062}
\cA_{m,n}^{\bth}\supseteq \bigcup_{\eta >
0}\La^{\bth}_{\TTT}(\psi^\eta\,) \,.
\end{equation}

\noindent Suppose $ \xmn\in\La^{\bth}_{\TTT}(\p^{\eta})$ for some
$\eta
> 0$.   By definition, (\ref{e:061}) is
satisfied for infinitely many $\ttt\in\TTT$. For any such $\ttt$,
there exists $\alpha=(\vv p,\vv q)\in \AAA$ such that
$$
    | g_{\vv t}M_{\xmn}^{\bth}\alpha|< 2^{-\eta \s t}\,.
$$
On taking the product of the first $m$ coordinates of
$g_{\ttt}M_{\xmn}^{\bth}\alpha$,  we obtain  that
$$
\prod_{j=1}^m2^{t_j}|\xmn_{j}\vv q+p_j+\theta_j| \, < \, 2^{-m\eta
\, \s t}\,.
$$
Similarly, the product of the last $n$ non-zero coordinates of
$g_{\ttt}M_{\xmn}^{\bth}\alpha$ yields that
$$
\prod_{\stackrel{\scriptstyle 1\le i\le n}{q_i\not=0}}2^{-t_{m+i}}
\;  |q_i| \, < \, 2^{-n \, \eta \,\s t}\,.
$$

\noindent  By  definition,  $ t_{m+i} \ \ge \, 0 $ $ (1\le i\le n)$
for  any $ \ttt = (t_1,\dots,t_{m+n}) \in \TTT $. Also, in view of
(\ref{e:052}) we have that $\s t \ge 0 $. Hence, by (\ref{e:039})
the above displayed  inequalities imply that
\begin{equation}\label{e:063}
\PI\big(\xmn\vv q+\vv p+\bth\big)<2^{-m\eta \,\s t-\s t/2} \qqand
\PI_+(\vv q)<2^{\s t/2}\,.
\end{equation}

\noindent By (\ref{e:053}), we have that (\ref{e:063}) is
satisfied for arbitrarily large $ \s t$. Now set $Q=2^{\s t/2}$
and $\ve:=2m\eta$. It follows that (\ref{e:054}) is satisfied for
arbitrarily large $Q$. The upshot is that
$\xmn\in\cA^{\bth}_{m,n}$.  This establishes (\ref{e:062}).

\smallskip

Steps 1 and 2 establish (\ref{e:046}) and  complete the proof of
Proposition \ref{prop5}.

\hfill $\boxtimes$

\vspace{1ex}

\noindent{\em Remark. } With
$\cA_{m,n}^{\bth}:=\{\xmn\in\Rmn:\wxth>1\}  $, as is the case when dealing with
Theorem~\ref{t1}A, the proof of Proposition \ref{prop5} remains
pretty much unchanged. The main difference is the manner in which we
define the set $\TTT$. Given $ s \in \Z_+  $ and $ l \in \Z_+ $, let
$\vv s :=(s,\dots,s) \in\Z_+^m$ and $\vv l :=(l\dots,l) \in\Z_+^n$.
On keeping the same notation as above, we have that $  \zeta=\frac{m
s- n l}{m+n}$. Furthermore, define $\vv t=(t_1,\dots,t_{m+n})$ by
setting
\begin{equation*}
    \vv t:=\big(\underbrace{s-\zeta,\dots,\ s-\zeta}_{m \text{ times}} \, , \ \underbrace{ l+\zeta,
    \dots,l+\zeta }_{n \text{ times}} \big)\,\tag{$\ref{e:048}\rule{1pt}{0pt}'$}
\end{equation*}
\noindent  and let
\begin{equation*}
    \TTT:=\{\ttt\in\R^{m+n}\text{ defined by $(\ref{e:048}\rule{1pt}{0pt}')$} \, : \,
    s\in\Z_+,\;
    l\in\Z_+  \text{ with } s \ge  t     \,     \}\, . \tag{$\ref{e:049}'$}
\end{equation*}

\noindent Note that $\TTT$  is a subset of the set defined by
(\ref{e:049}). Thus, conditions (\ref{e:039}) and (\ref{e:045}) are
automatically satisfied for this `smaller' choice of $\TTT$. To
establish the set equality (\ref{e:046}), we start by verifying that
$\xmn\in\cA^{\bth}_{m,n}$ if and only if there exists an $\ve>0$,
such that for arbitrarily large $Q>1$ there is an $\alpha=(\vv p,\vv
q)\in\cA=\Z^m\times(\Z^n\bnz)$ satisfying $
 \|\xmn\vv q+\vv p+\bth\|^m<Q^{-1-\ve} $ and  $ |\vv q|^n\le Q $.
These inequalities replace those appearing in (\ref{e:054}) and by
naturally modifying the arguments setout in Steps 1 and 2 above, we
obtain (\ref{e:046}).

\section{An Inhomogeneous Transference Principle}\label{HIT}

In this section we  develop a general framework that allows us to
transfer zero measure statements for homogeneous $\limsup$ sets to
inhomogeneous $\limsup$ sets.  To a certain extent, the framework is motivated
by our desire to establish the specific transference given by (\ref{e:047})
 and thereby complete the proof of Theorem \ref{t1}.

\medskip

Let $(\Omega,d)$ be a locally compact metric space. Given two
countable `indexing' sets $\AAA$ and $\TTT$, let $\hom$ and $\inh$
be two maps from $\TTT\times \AAA\times\Rp $ into the set of open
subsets of $\Omega$ such that
$$
\hom\,:\,(\ttt,\alpha,\ve)\in \TTT\times \AAA\times\Rp \,\mapsto
\,\hom_\ttt(\alpha,\ve)
$$
and
$$
\inh\,:\,(\ttt,\alpha,\ve)\in \TTT\times \AAA\times\Rp \,\mapsto
\,\inh_\ttt(\alpha,\ve)\,.
$$

\noindent Furthermore, let

\begin{equation}\label{e:064}
\hom_\ttt(\ve):=\bigcup_{\alpha\in \AAA}\hom_\ttt(\alpha,\ve)\qqand
\inh_\ttt(\ve):=\bigcup_{\alpha\in \AAA}\inh_\ttt(\alpha,\ve)\,.
\end{equation}

\noindent Next, let $\bm\Psi$ denote a set of functions $
\psi:\TTT\to\Rp\,:\,\ttt\mapsto \psi_\ttt\,. $ For $\p\in\bm\Psi$,
consider the $\limsup$ sets

\begin{equation}\label{e:065}
 \La_\hom(\psi\,)=\limsup_{\ttt \in \TTT}\hom_\ttt(\p_\ttt)
 \qqand
 \La_\inh(\psi\,)=\limsup_{\ttt \in \TTT}\inh_\ttt(\p_\ttt)\,.
\end{equation}
\bigskip

\noindent For reasons that will soon become apparent, we refer to
sets associated with the map $\hom$  as homogeneous sets and those
associated with the map $\inh$ as inhomogeneous sets. The following
`intersection' property states that the intersection of two distinct
inhomogeneous sets is contained in a homogeneous set.

\bigskip

\noindent\textbf{The intersection property. } The triple
$(\hom,\inh,\bm\Psi)$ is said to satisfy \emph{the intersection
property} if for any $\psi\in\bm\Psi$, there exists
$\psi^*\in\bm\Psi$ such that for all but finitely many $\ttt\in\TTT$
and all  distinct $\alpha$ and $\alpha'$ in $\AAA$ we have that
\begin{equation}\label{e:066}
    \inh_\ttt(\alpha,\p_\ttt)\cap \inh_\ttt(\alpha',\p_\ttt)\subset
\hom_{\ttt}(\p^*_{\ttt})   \ .
\end{equation}

\bigskip

\noindent The following notion of `contracting' is the natural
generalisation of  the $\Rmn$ version stated in  \S\ref{contr}.

\bigskip

\noindent\textbf{The contracting property. } Let $\mu $ be a
non-atomic, finite, doubling  measure supported on a bounded subset
$\Supp$ of $\Omega$.  We say that  $\mu$ is \emph{contracting with
respect to\/ $(\,\inh,\bm\Psi)$}\/ if for any $\psi\in\bm\Psi$ there
exists $\psi^+\in\bm\Psi$ and a sequence of positive numbers
$\{k_\ttt\}_{\ttt\in \TTT}$ satisfying
\begin{equation}\label{e:067}
    \sum_{\ttt\in \TTT}k_\ttt<\infty  \ ,
\end{equation}
 such that for all but finitely  $\ttt\in \TTT$ and all
$\alpha\in \AAA$ there exists a collection $\Cta$ of balls $B$
centred at $\Supp$ satisfying the
following conditions\,{\rm:}
  \begin{equation}\label{e:068}
    \Supp\cap\inh_\ttt(\alpha,\psi_\ttt) \ \subset \
    \bigcup_{B\in\Cta}B\
  \end{equation}
    \begin{equation}\label{e:069}
        \Supp\cap\bigcup_{B\in\Cta}B \ \subset \ \inh_\ttt(\alpha,\p^+_{\ttt})
    \end{equation}
    and
    \begin{equation}\label{e:070}
        \mu\Big(5B\cap\inh_\ttt(\alpha,\psi_\ttt)\Big)\ \le  \ k_\ttt\ \,  \mu(5B) \ .
    \end{equation}

\vspace{4ex}

The intersection and contracting properties enable us to transfer
zero $\mu$-measure statements for the homogeneous $\limsup$ sets
$\La_\hom(\psi\,)$ to the inhomogeneous $\limsup$ sets
$\La_\inh(\psi\,) $.

\begin{theorem}\label{t5} {\bf (Inhomogeneous Transference
Principle)} \  Suppose\/ that $(\hom,\inh,\bm\Psi)$ satisfies the
intersection property and that $\mu$ is contracting with respect
to\/ $(\inh,\bm\Psi)$.  Then
\begin{equation}\label{e:071}
\mu(\La_\hom(\psi))=0    \ \quad \forall\ \psi\in\bm\Psi\
\qquad\Longrightarrow\qquad
 \mu(\La_\inh(\psi))=0 \ \quad \forall\ \psi\in\bm\Psi  \,.
\end{equation}
\end{theorem}

\bigskip

Before proving Theorem \ref{t5}, we consider an
application that establishes (\ref{e:047}) and at the same time
 clarifies the above abstract setup.

\subsection{Completing the proof of Theorem~\ref{t1}}\label{me}

Given the  Inhomogeneous Transference Principle, we are fully
armed to complete the proof of Theorem~\ref{t1}. In view of the
reformulation and reduction carried out in \S\ref{sec:1}, both
parts of  Theorem \ref{t1} follow on establishing (\ref{e:047})
with an appropriate choice of $\TTT$. As in \S\ref{sec:1}, we will
concentrate on the proof of Theorem \ref{t1}B -- part (B) of
Theorem \ref{t1}.

Throughout $\bm\theta\in\R^m$ is fixed.  Let $\mu$ be a measure on
$\Rmn$ that is strongly contracting almost everywhere and fix a
set $\TTT$ arising from Proposition \ref{prop5}. In terms of
establishing (\ref{e:047}), sets of $\mu$-measure zero are
irrelevant. Therefore we can simply assume that $\mu$ is strongly
contracting. We  show that (\ref{e:047}) falls within the scope of
the above general framework. Let $\Omega := \Rmn$ and let $\cA$ be
given by (\ref{e:041}). Given $\ve \in \Rp$,  $\ttt \in \TTT$ and
$ \alpha \in \cA$  let
$$
\hom_\ttt(\alpha,\ve):=\De_\ttt(\alpha,\ve)=\De^{\vv0}_\ttt(\alpha,\ve)\qquad\text{and}\qquad
\inh_\ttt(\alpha,\ve):=\De^{\bth}_\ttt(\alpha,\ve),
$$
where $\De^{\bth}_\ttt(\alpha,\ve)$ is defined by (\ref{e:042}).
This defines the  maps $\hom$ and $\inh$ associated with the
general framework.  It is readily seen that
$\hom_\ttt(\ve)=\De^{\vv0}_\ttt(\ve)$ and
$\inh_\ttt(\ve)=\De^{\bth}_\ttt(\ve)$. Next, let $\bm\Psi$ be the
class of  functions given by (\ref{e:043}). Then,  it immediately
follows that
$$
  \La_\hom(\psi) \,  = \, \La_{\TTT}(\p) := \La^{\vv0}_{\TTT}(\p)
  \qquad\text{and}\qquad \La_\inh(\psi) \, = \,   \La^{\bth}_{\TTT}(\p) \ ,
$$
where the set $\La^{\bth}_{\TTT}(\p)$ is  defined by
(\ref{e:044}). The upshot is that (\ref{e:047}) and  (\ref{e:071})
are precisely the same statement. In view of the Inhomogeneous
Transference Principle,  it follows that (\ref{e:047})  is a
consequence of verifying that $(\hom,\inh,\bm\Psi)$ satisfies the
intersection property and that $\mu$ is contracting with respect
to $(\inh,\bm\Psi)$.

\bigskip

\noindent\text{\em Verifying the intersection property. } Let $\psi\in\bm\Psi$.
This means that $\psi_\ttt=2^{-\eta\,\s t}$ for some constant
$\eta>0$. To establish the intersection property given by
(\ref{e:066}),  define $\p^*$ by setting $\psi^*_\ttt=2^{-\eta/2\,\s
t}$. Obviously,  $\psi^*\in\bm\Psi$. Next, fix two distinct
 $\alpha$ and $\alpha'$ in $\cA$.  By definition,
$$
\alpha=(\vv p,\vv q)\qquad\text{and}\qquad\alpha'=(\vv p',\vv q')
$$
for some $\vv p,\vv p'\in\Z^m$ and $\vv q,\vv q'\in\Z^n\bnz$. Now, take
any point
$$\xmn \ \in \ \inh_\ttt(\alpha,\p_\ttt)\cap
\inh_\ttt(\alpha',\p_\ttt)\, := \,  \De^{\bth}_\ttt(\alpha,\p_\ttt)  \cap    \De^{\bth}_\ttt(\alpha',\p_\ttt)   \ . $$
Clearly, we may as well assume that the intersection is non-empty. In the notation of
\S\,\ref{sec:1}, we have that
\begin{equation}\label{e:072}
| g_{\vv t}M_{\vv X}^{\bth}\alpha|\,< \p_\ttt\qquad\text{and}\qquad
| g_{\vv t}M_{\vv X}^{\bth}\alpha'|\,< \p_\ttt\,.
\end{equation}
Let $\alpha'':=(\vv p'',\vv q'')$, where   $\vv p'':=\vv p-\vv p'\in\Z^m$
and $\vv q'':=\vv q-\vv q'\in\Z^n$. Since $\TTT$ satisfies  (\ref{e:045}) and therefore (\ref{e:053}), it follows that
\begin{equation}\label{e:073}
| g_{\vv t}M_{\vv X}\alpha''| \, = \,  | g_{\vv t}M_{\vv
X}(\alpha-\alpha')| \, = \, |g_{\vv t}M_{\vv X}^{\bth}\alpha  -g_{\vv
t}M_{\vv X}^{\bth}\alpha'|\,\stackrel{(\ref{e:072})}{<} \,
2\p_\ttt  \, <  \, \p^*_\ttt
\end{equation}
for all but finitely many $\ttt\in\TTT$.
If $\vv q''=\vv0$,  we obtain via
(\ref{e:073})  that $|\vv p''|<\p^*_\ttt<1$ for
all but finitely many $\ttt\in\TTT$.  However $\vv p''\in\Z^m$ and so we must have
$\vv p''=\vv0$. This contradicts the assumption that
$\alpha\not=\alpha'$. The upshot is that $\vv q''\not=\vv0$ and
so $\alpha''\in\cA$. Hence, it follows that
$$\xmn\in
\De_\ttt(\alpha'',\p^*_\ttt)\subset \De_\ttt(\p^*_\ttt)  := \hom_\ttt(\p^*_\ttt)$$
for all but finitely many $\ttt\in\TTT$. This verifies the intersection property.

\bigskip

\noindent\text{\em Verifying the contracting property. } Recall,
$\mu$ is a measure on $\Rmn$ that is strongly contracting.
Therefore, $\mu$ is by definition  non-atomic, doubling and
finite. Also without loss of generality we can assume that the
support $\Supp$ of $\mu$ is bounded. Thus, to establish that $\mu$
is contracting with respect to $(\inh,\bm\Psi)$ it remains to
verify the conditions given by  (\ref{e:067}) -- (\ref{e:070}).
Fix  $\p\in\bm\Psi$. Then,  $\psi_\ttt=2^{-\eta\,\s t}$ for some
constant $\eta>0$ and we define $\p^+$ by setting $\p^+_\ttt
:=\sqrt {\p_\ttt}$. Obviously,  $\psi^+\in\bm\Psi$.  Let $r_0$ be
the positive constant appearing in the definition of strongly
contracting. Since $\TTT$ satisfies (\ref{e:045}) and therefore
(\ref{e:053}), it follows that
\begin{equation}\label{e:074}
\p^+_\ttt \le   \min \{ 1, r_0 \}  \qquad\text{and}\qquad \s t \ge 0
\end{equation}
for all but finitely many $\ttt \in \TTT$. Now fix such a
$\ttt=(t_1,\dots,t_{m+n})\in\TTT$ and $\alpha'=(\vv p,\vv
q)\in\cA$.  The set $\inh_\ttt(\alpha',\p_\ttt)$ corresponds to
$\xmn \in \Rmn$ satisfying
\begin{equation}\label{e:075}
|\xmn_j \, \vv q+p_j+\theta_j|<2^{-t_j}\p_\ttt\quad(1\le j\le
m)\quad\text{with}\quad |q_i|<2^{t_{m+i}}\p_\ttt\quad(1\le i\le n)
\ .
\end{equation}
Similarly, $\inh_\ttt(\alpha',\p^+_\ttt)$ corresponds  to $\xmn
\in \Rmn$ satisfying
\begin{equation}\label{e:076}
|\xmn_j \, \vv q+p_j+\theta_j|<2^{-t_j}\p^+_\ttt\quad(1\le j\le
m)\quad\text{with}\quad |q_i|<2^{t_{m+i}}\p^+_\ttt\quad(1\le i\le
n)  \ .
\end{equation}

\noindent Without loss of generality, we assume that the right
hand side inequalities of (\ref{e:075}) and (\ref{e:076}) are
fulfilled for $\alpha'=(\vv p,\vv q)$.  Otherwise, the sets under
consideration are empty and the conditions (\ref{e:068}) --
(\ref{e:070}) are easily met. For $  j \in \{ 1, \ldots ,  m \} $,
define
\begin{equation}\label{e:077}
 \ve_j=\ve_{j,\ttt}:=2^{-t_{j}}\p^+_\ttt/|\vv q|_2\qquad   \text{and \ let }  \qquad   \delta=\delta_\ttt:=\p^+_\ttt.
\end{equation}

\noindent By (\ref{e:039}) and the fact that $\s t \ge  0 $, it
follows  that $\sum_{j=1}^mt_j\ge0$ and so there exists
$j\in\{1,\dots,m\}$ such that $2^{-t_j}\le1$. Since $\vv
q\in\Z^n\bnz$ we have that $|\vv q|_2^{-1}\le 1$. Therefore
$\min_{1\le j\le m}\ve_{j,\ttt}<\p^+_\ttt$ and  (\ref{e:074})
implies that
$$
 \min_{1\le j\le m}\ve_{j,\ttt}< r_0\qqand \delta_\ttt< 1   \ .
$$

\noindent By the definition of $\p^+_\ttt$, we have that
$\delta\ve_j=2^{-t_{j}}\p_\ttt/|\vv q|_2$ for each $  j \in \{ 1,
\ldots ,  m \} $. Therefore, by (\ref{e:075}) and (\ref{e:076}),
it follows that
\begin{equation}\label{e:078}
\inh_\ttt(\alpha',\p_\ttt)=\cL_{\vv a,\vv
b}^{(\delta\bm\ve)}\qquad\text{and}\qquad
\inh_\ttt(\alpha',\p^+_\ttt)=\cL_{\vv a,\vv b}^{(\bm\ve)}   \ ,
\end{equation}
 where
$\cL_{\vv a,\vv b}^{(\delta\bm\ve)}$ and $\cL_{\vv a,\vv
b}^{(\bm\ve)}$ are defined by (\ref{e:012}) with $\vv a:=\vv
q/|\vv q|_2$ and $\vv b:=(\vv p+\bm\theta)/|\vv q|_2$.
 By the definition of strongly contracting, for each
$\xmn\in\cL^{(\delta\bm\ve)}_{\vv a,\vv b}\cap\Supp$ there is an
open  ball $B=B_{\xmn}$ centred at $\xmn$ such that (\ref{e:013})
and (\ref{e:014}) are satisfied. With reference to the general
framework,  define ${\cal C}_{\ttt, \alpha'}$ to be the collection
of all such balls. By definition, each point $\xmn\in
\cL^{(\delta\bm\ve)}_{\vv a,\vv b}\cap\Supp$ is the centre of a
ball in ${\cal C}_{\ttt, \alpha'}$. Finally, for any $\ttt \in
\TTT$,  let
$$k_\ttt \; := \; C \, (\p^+_\ttt)^\alpha   \  $$ where $C $ and $ \alpha$ are the
constants appearing in  the definition of strongly contracting.
Then

$$
\begin{array}{l}
 \text{(\ref{e:067}) follows from (\ref{e:045});}
\\[1ex]
 \text{(\ref{e:068}) follows from the definition of ${\cal C}_{\ttt, \alpha'}$;}
\\[1ex]
 \text{(\ref{e:069}) and (\ref{e:070}) follow from (\ref{e:013}) and (\ref{e:014}) via (\ref{e:078}).}
\end{array}
$$

\smallskip

\noindent This verifies that $\mu$ is contracting with respect to
$(\inh,\bm\Psi)$.

\hfill $\boxtimes$

\bigskip

\noindent {\it Remark.} When dealing with Theorem~\ref{t1}A, the
above arguments remain essentially unchanged. The main difference
is that we work with the specific  $\TTT$  given by
$(\ref{e:049}\,')$. Then for any $\ttt
=(t_1,\dots,t_{m+n})\in\TTT$, we have that  $t_1=\dots=t_m$. This
 ensure that $\ve_1,\dots,\ve_m$ as defined by
(\ref{e:077}) are all equal and therefore $\mu$ being contracting
rather than strongly contracting is sufficient.

\section{Preliminaries for Theorem \ref{t5}}

In this section we group together various self contained
statements  that we appeal to during the course of establishing
the Inhomogeneous Transference Principle.  We start with a basic
covering result from geometric measure theory.

\begin{lemma} \label{lemma1}
Every collection  $\cC$ of balls of uniformly bounded diameter in
a metric space  contains a disjoint subcollection $\cG$ such that
$$
\bigcup_{B\in\cC}B\subset\bigcup_{B\in\cG}5B\,.
$$
\end{lemma}

 This covering lemma is usually referred to as the $5r$-lemma.
For further details and proof the reader is refereed to
\cite{Heinonen,Mattila-1995}. The following  enables us to bound the cardinality of the disjoint collection
arising from the $5r$-lemma.

\begin{lemma}\label{mewlemma}
Let $(\Omega,d)$ be a metric space equipped with a finite measure
$\mu$. Then every disjoint collection $\cC$ of $\mu$-measurable
subsets of $\Omega$ with positive $\mu$-measure is at most countable.
\end{lemma}

\smallskip

\noindent{\em Proof. \ } For $k\in\Z$, let
$\cC^{(k)}$ be the subcollection of $\cC$ consisting of sets $B\in\cC$
such that $2^{k}\le\mu(B)<2^{k+1}$. Obviously we have that $\cC=\bigcup_{k\in\Z}\cC^{(k)}$.
Since the balls in $\cC$ are pairwise disjoint, it follows  that
$\#\cC^{(k)}\le2^{-k}\mu(\Omega)<\infty$. Therefore, $\cC$ is a countable union of finite sets and so is at most countable.

\hfill $\boxtimes$

 The next statement is a simple
consequence of the continuity of measures.

\begin{lemma}\label{lemma2}
Let $(\Omega,d)$ be a metric space equipped with a measure $\mu$.
 Let $\{A_i\}_{i\in\N}$ be a sequence of
$\mu$-measurable subsets of $\Omega$ and
$$A_\infty \; :=  \; \limsup_{i\to\infty}A_i=\bigcap_{m=1}^\infty\bigcup_{i=m}^\infty
A_i   \ . $$ Then $\mu(A_\infty)=0$ if and only if for any
$\ve>0$, there exists a positive constant $m_0(\ve)$ such that
$\mu(\bigcup_{i=m}^\infty A_i)<\ve$ for all  $m>m_0(\ve)$.
\end{lemma}

\noindent Recall, that the Borel-Cantelli lemma from probability
theory states that
$$
\mu(A_\infty)=0   \qquad \text{if} \qquad \sum_{i=1}^\infty \mu
(A_i) \, < \, \infty \ .
$$

Sprindzuk's proof of Mahler's conjecture is based on the notions
of essential and inessential domains -- see \cite[\S14]
{Sprindzuk-1979-Metrical-theory}. As we shall soon see, the proof
of Theorem \ref{t5} is based on the related notions of
$\mu$-essential and $\mu$-inessential  balls. The
following lemma enables us to exploit these key notions.

\begin{lemma}\label{lemma3}
Let $(\Omega,d)$ be a metric space equipped with a finite measure
$\mu$. Let $\{B_i\}_{i\in\N}$ be a sequence of $\mu$-measurable
subsets of\/ $\Omega$. Suppose there exists a constant $c\in(0,1)$
such that for every $i\in\N$
\begin{equation}\label{e:079}
\textstyle\mu\Big(B_i\cap\bigcup_{j\not=i}B_j\Big)\le
c\,\mu(B_i)\,.
\end{equation}

\noindent Then
$$
\sum_{i=1}^\infty\mu(B_i)\ \le \ \frac1{1-c}\ \mu(\Omega)\,.
$$
\end{lemma}

\medskip

\noindent{\em Proof. \ } Given $i\in\N$, let
$$
B_i^{\rm 0}:=B\setminus\bigcup_{j\not=i}B_j\qquad\text{and}\qquad
B_i^{\rm 1}:=B\cap\bigcup_{j\not=i}B_j\,.
$$
Thus, $B_i^{\rm 0}$ corresponds to the region of $B_i$ that is
disjoint from the  sets  $B_j$ with $j\not=i$.  On the other hand,
$B_i^{\rm 1}$ corresponds to the region of $B_i$ that is
non-disjoint from the sets  $B_j$ with $j\not=i$.   Obviously
$B_i^{\rm 0}$ and $B_i^{\rm 1}$ are $\mu$-measurable. By
(\ref{e:079}), we have that $\mu(B_i^0)=\mu(B)-\mu(B_i^1) >
(1-c)\mu(B)$. Thus,
\begin{equation}\label{e:080}
    \mu(B)\le\frac1{1-c}\ \mu(B_i^0)\,.
\end{equation}
By construction,  $B_i^0\cap B_j^0=\emptyset$ for distinct
$i,j\in\N$. Therefore
$$
\sum_{i=1}^\infty\mu(B_i)\stackrel{(\ref{e:080})}{\le}\frac{1}{1-c}\
\sum_{i=1}^\infty\mu(B_i^0) \ = \ \frac{1}{1-c}\
\mu\Big(\bigcup^\circ_{i\in\N}B_i^0\Big) \ {\le}\ \frac1{1-c}\
\mu(\Omega)\,.
$$
\hfill $\boxtimes$

\section{Proof of Theorem~\ref{t5}}\label{proof_t5}

Fix any $\psi \in \Psi$. The goal is to show that
\begin{equation}\label{e:081}
\mu (\La_\inh(\psi) )   = 0   \, .
\end{equation}
We are given that $\mu$ is contracting with respect to
$(\inh,\bm\Psi)$. Note that we can assume that this contacting
property and indeed the intersection property defined in
\S\ref{HIT},   are valid for all $\ttt \in \TTT $ rather than all
but finitely many  --  removing a finite number of elements from
$\TTT$ does not alter the $\limsup$ set $\La_\inh(\psi)$ under
consideration. With this in mind, let $\psi^+\in\bm\Psi$ be the
function, $\{k_\ttt\}$ be the sequence  and $\Cta$ be the
collection of balls arising  from  the contracting property. Since
$\Supp=\supp\mu$ is bounded, we can assume that the balls in
$\Cta$ are of radius bounded by $r=\diam(\Supp)$. Indeed, if
$\Cta$ contains a ball $B$ of radius bigger than $r$ centered at
$x$ then we can replace $B$ with $B(x,r)$. This replacement would
not affect the properties (\ref{e:068})--(\ref{e:070}) as
$\Supp\subset B(x,r)$. Thus, in view of Lemma~\ref{lemma1} there
is a disjoint subcollection $\Gta$ of $\Cta$ such that
\begin{equation}\label{e:082}
    \bigcup_{B\in\Cta}B\subset \bigcup_{B\in\Gta}5B\,.
\end{equation}
Every ball $B\in\Gta$ is centred at the support of $\mu$ and so is of positive
$\mu$-measure. Thus, by Lemma~\ref{mewlemma}, the collection $\Gta$ is at
most countable.  In view of (\ref{e:068}) and (\ref{e:082}) we have that
$$
\Supp\cap\iDta\ \subset \ \bigcup_{B\in\Gta}5B\cap \iDta\ \subset \
\iDta\,.
  $$
Taking the union over all $\alpha\in\cA$ and using (\ref{e:064}) gives
$$
    \Supp\cap
    \inh_\ttt(\p_\ttt)\
    \subset \ \bigcup_{\alpha\in \AAA}\bigcup_{B\in\Gta}5B\cap \iDta\ \subset \ \inh_\ttt(\p_\ttt)\,.
$$
In view of (\ref{e:065}), it follows that
$$
    \Supp\cap\La_\inh(\psi)\ \subset \ \limsup_{\ttt \in \TTT}\bigcup_{\alpha\in
\AAA}\bigcup_{B\in\Gta}5B\cap \iDta\ \subset\ \La_\inh(\psi)\,.
$$
This implies that
\begin{equation}\label{e:083}
    \mu(\La_\inh(\psi))\ =\ \mu\Big(\limsup_{\ttt \in \TTT}\bigcup_{\alpha\in
\AAA}\bigcup_{B\in\Gta}5B\cap \iDta\Big)\,.
\end{equation}
Thus, (\ref{e:081}) will follow on establishing that the right hand
side of (\ref{e:083}) is zero. The key that enables us to do
precisely this, is the following decomposition of the right hand
side $\limsup$ set in terms of $\mu$-essential and
$\mu$-inessential balls.

\medskip

\noindent{\bf Definition. } A ball  $B\in \Gta$ is said to be
\emph{$\mu$-essential} if
\begin{equation}\label{e:084}
\mu\Big(B\
\cap\!\!\bigcup_{\alpha'\in\cA\smallsetminus\{\alpha\}}\bigcup_{B'\in\cG_{\ttt,\alpha'}}B'\Big)\
\le \ \frac 12 \ \mu(B)
\end{equation}
 and \emph{$\mu$-inessential}
otherwise.

\medskip

\noindent Let $\Dta$  denote  the collection of $\mu$-essential
balls in $\Gta$ and let $\Nta$ denote the collection of
$\mu$-inessential balls in $\Gta$. Consider the corresponding
limsup sets
$$
\Lad:=\limsup_{\ttt \in \TTT} \ \ \bigcup_{\alpha\in \AAA}\ \
\bigcup_{B\in\Dta}5B\cap \iDta\, ,
$$
and
$$
\Lan=\limsup_{\ttt \in \TTT }\ \ \bigcup_{\alpha\in \AAA}\ \
\bigcup_{B\in\Nta}5B\cap \iDta\,.
$$
It is easily  seen that the limsup set in the right hand side of
(\ref{e:083}) is equal to  $\Lan\cup\Lad$ and so we have that
$$
\mu(\La_\inh(\psi))\ =\mu(\Lad\cup\Lan)  \ .
$$

\noindent Thus, the  statement of Theorem~\ref{t5} will follow on showing that
$\mu(\Lad)=0$ and $\mu(\Lan)=0$.

\subsection{The $\mu$-essential case: $\mu(\Lad)=0$}

In view of (\ref{e:084}), Lemma~\ref{lemma3} and the fact that $\Gta$ is a disjoint
collection of balls,  we have that
\begin{equation}\label{e:085}
\sum_{\alpha\in\cA}\sum_{B\in\Dta}\mu(B)\le2\mu(\Omega)\,.
\end{equation}

\noindent It follows that

$$
\begin{array}[b]{rcl}
  \displaystyle\mu\Big(\bigcup_{\alpha\in \AAA}
 \ \bigcup_{B\in\Dta}5B\cap \iDta\Big) & \le & \displaystyle\sum_{\alpha\in
 \AAA}\ \sum_{B\in\Dta}\mu\Big(5B\cap \iDta\Big) \\[4ex]
   & \stackrel{(\ref{e:070})}{\le} & \displaystyle\sum_{\alpha\in \AAA}\sum_{B\in\Dta}k_\ttt\ \mu(5B) \\[4ex]
   & \stackrel{(\ref{e:010})}{\le} & k_\ttt\ \lambda^3\displaystyle\sum_{\alpha\in
 \AAA}\ \sum_{B\in\Dta}\mu(B)
   \\[4ex]
   & \stackrel{(\ref{e:085})}{\le} &\displaystyle2\ \mu(\Omega)\lambda^3k_\ttt\, .
\end{array}
$$

\noindent Since $\mu(\Omega)<\infty$, the quantity  $2
\mu(\Omega)\lambda^3$ is a finite positive constant. Hence,
$$
\sum_{\ttt\in \TTT}  \ \; \mu\Big(\bigcup_{\alpha\in
\AAA}\bigcup_{B\in\Dta}5B\cap \iDta\Big)  \  \ll  \
  \sum_{\ttt \in \TTT}k_\ttt  \,  \stackrel{(\ref{e:067})}{<}   \,   \infty\,.
$$

\noindent A simple consequence of the Borel-Cantelli lemma  is that
$
\mu(\Lad)=0\,.
$

\subsection{The $\mu$-inessential  case: $\mu(\Lan)=0$}

For a ball $B\in\Nta$, let
$$
\Bn:=B\
\cap\!\!\bigcup_{\alpha'\in\cA\smallsetminus\{\alpha\}}\bigcup_{B'\in\cG_{\ttt,\alpha'}}B'\,.
$$
By definition,
\begin{equation}\label{e:086}
\mu(\Bn)>\tfrac12\,\mu(B)   \ .
\end{equation}
Consider the following two sets
$$
\Lan'= \limsup_{\ttt \in \TTT
}\bigcup_{\alpha\in\cA}\bigcup_{B\in\Nta}5B\qqand \Lan''=
\limsup_{\ttt \in \TTT
}\bigcup_{\alpha\in\cA}\bigcup_{B\in\Nta}\Bn\,.
$$
Obviously we have that
$
\Lan\subset \Lan'\, .
$
Hence, it suffices to show that
$$\mu(\Lan')=0  \, . $$

\medskip

\noindent \textbf{Step 1:  $\mu(\Lan'')=0$. } Note that in view of
(\ref{e:069}),  for any $B\in \Nta \subset \Gta$ we have that
\begin{equation}\label{e:087}
        \Supp\cap B \; \subset \; \inh(\alpha,\p^+_\ttt) \ .
\end{equation}
By definition,  if $x\in\Supp\cap\Bn$ then there exists
$\alpha'\in \AAA\smallsetminus\{\alpha\}$ and a ball
$B'\in\cG_{\ttt,\alpha'}$ such that $x\in B'$. Again, by
(\ref{e:069}) we have that
\begin{equation}\label{e:088}
        \Supp\cap B' \; \subset \; \inh(\alpha',\p^+_\ttt)  \ .
\end{equation}
Since $x\in \Supp\cap B\cap B'$, it follows via (\ref{e:087}) and
(\ref{e:088}) that
$$x\in \inh(\alpha,\p^+_\ttt)\cap\inh(\alpha',\p^+_\ttt).$$
We are given that $(\hom,\inh,\bm\Psi)$ satisfies the intersection
property. Thus, in view of (\ref{e:066}) we conclude that
$$x\in\hom_{\ttt}(\psi^*_{\ttt})  \ , $$  where $\p^*\in\bm\Psi$ is
associated with $\p^+$.   The upshot  is that if $ x \in \Supp\cap
\Lan''$, then $x$ lies in the `homogeneous' sets
$\hom_{\ttt}(\psi^*_{\ttt})$ for infinitely many $\ttt \in \TTT$.
In other words, $\Supp\cap \Lan''  \subset  \La_\hom(\p^*)$.
However, homogeneous $\limsup$ sets are assumed to be of
$\mu$-measure zero - see (\ref{e:071}). Therefore,
\begin{equation}\label{e:089}
       \mu(\Lan'')=0  \ .
\end{equation}

\bigskip

\noindent \textbf{Step 2: $\mu(\Lan')=0$. } Fix any enumeration
$\{\ttt_l\}_{l\in\N}$ of the set $\TTT$. For $m\in\N$, let
\begin{equation}\label{e:090}
\Lan''(m):=\bigcup_{l\ge m}\ \bigcup_{\alpha\in\cA} \bigcup_{B\in
\cN_{\ttt_l,\alpha}}\Bn\,.
\end{equation}
 Then,  $\Lan''=\bigcap_{m=1}^\infty
\Lan''(m)$.  In view of (\ref{e:089}), it follows  via
Lemma~\ref{lemma2} that for any $\delta>0$, there exists a
positive constant $m_0(\delta)$ such that
\begin{equation}\label{e:091}
\mu\big(\Lan''(m)\big)<\delta\  \qquad \forall \quad m>m_0(\delta)
\ .
\end{equation}
Let $\cG(m)$ be the collection of balls $5B$ such that
$B\in\cN_{\ttt_l,\alpha}$ with $l\ge m$ and $\alpha\in\cA$. By
(\ref{e:090}), for any $5B\in\cG(m)$ we have that
$B\cap\Lan''(m)\supset \Bn$. Therefore,
\begin{equation}\label{e:092}
    \mu(B\cap\Lan''(m))\ge  \mu(\Bn)   \stackrel{(\ref{e:086})}{>}
    \tfrac12\mu(B)
\end{equation}
for any $5B\in\cG(m)$. Recall that  $\Supp$ is bounded and so
$\cG(m)$ is a collection of balls of uniformly bounded diameter.
Thus, in view of  Lemma~\ref{lemma1} there is a disjoint
subcollection $\cG'(m)$ of $\cG(m)$ such that
$$
\bigcup_{5B\in\cG(m)}5B \ \ \subset\bigcup_{5B\in\cG'(m)}25B\, .
$$
Finally, for $m\ge m_0(\delta)$  let
$$
\Lan'(m) \ := \ \bigcup_{l\ge m}\ \bigcup_{\alpha\in\cA} \
\bigcup_{B\in \cN_{\ttt_l,\alpha}}5B\,.
$$
Then
\begin{eqnarray*}
  \mu(\Lan'(m)) &\le& \mu\Big(\bigcup_{5B\in\cG'(m)}25B\Big)
  \\[2ex]
   &\le& \sum_{5B\in\cG'(m)}\mu\big(25B\big) \\
   &\stackrel{(\ref{e:010})}{\le}& \lambda^5 \sum_{5B\in\cG'(m)}\mu\big(B\big) \\
   &\stackrel{(\ref{e:092})}{\le}& 2\lambda^5 \sum_{5B\in\cG'(m)}\mu\big(B\cap\Lambda''_{\rm
   N}(m)\big)\\
   &\stackrel{}{=}& 2\lambda^5 \mu\Big(\bigcup_{5B\in\cG'(m)}^\circ B\cap\Lambda''_{\rm
   N}(m)\Big)\\
   &\stackrel{}{\le}& 2\lambda^5 \mu\big(\Lambda''_{\rm
   N}(m)\big)\\ \bigskip
   &\stackrel{(\ref{e:091})}{\le}& 2\lambda^5 \delta\,.
\end{eqnarray*}
Since $\Lan'=\bigcap_{m=1}^\infty \Lan'(m)$, it follows via
Lemma~\ref{lemma2} that $\mu(\Lan')=0$.

\hfill $\boxtimes$


\section{Final comments and open problems}

In principle, there are numerous problems that can be treated by
applying the basic recipes introduced in this paper -- in
particular the Inhomogeneous Transference Principle of
\S\ref{HIT}. The goal here is to indicate the diversity of these
problems beyond those considered in the main bulk of the paper.

\bigskip\noindent\textbf{Beyond simultaneous and dual extremality: $d$-extremality. }
The dual theory of Diophantine
approximation is concerned with approximation of points in $\R^n$ by
$(n-1)$--dimensional rational planes; i.e. rational hyperplanes. The
simultaneous theory of Diophantine approximation is concerned with
approximation of points in $\R^n$ by $0$--dimensional rational
planes; i.e. rational points. As a consequence of Khintchine's
transference principle, the two forms of approximation  lead to
equivalent notions of extremality in the homogeneous case. However,
as a consequence of a recent work by Laurent a lot more is true. For
$d \in \{0, \ldots n-1\}$, it is natural to consider the Diophantine
approximation theory in which points in $\R^n$ are approximated by
$d$-dimensional rational planes -- the dual and simultaneous
theories just represent the extreme. The related homogeneous
Diophantine exponents $\omega_d(\vv x)$ have been studied in some
depth by Schmidt \cite{Schmidt-67:MR0213301} in the sixties and more
recently by Laurent \cite{Laurent}. We refer the reader to Laurent
\cite{Laurent} for the definition of these exponents. For the
purpose of this discussion, it suffices to say that for any $d \in
\{0, \ldots n-1\}$, we have that $ \omega_d(\vv x) \ge (d+1)/(n-d)$
for all $\vv x \in\R^n.$ Now let $\mu$ be a measure on $\R^n$ and
say that $\mu$ is {\em $d$--extremal}  if
$$
\text{  $\omega_d(\vv x)=(d+1)/(n-d)$ \qquad for $\mu$-almost all
$\vv x \in \R^n$. }
$$
Khintchine's transference principle implies that the extreme cases
($d=0$ and $d=n-1$) of $d$--extremal are equivalent. Leaving aside
the details, the Schmidt-Laurent `Going-up' and `Going-down'
transference inequalities imply that all $n$ notions of
$d$--extremality  are in fact equivalent. This together with Theorem
KLW implies the following statement.

\begin{theorem} \label{t10}
Any friendly measure on $\R^n$ is $d$-extremal for all $d \in \{0,
\ldots n-1\}$ .
\end{theorem}

\noindent Given the equivalence of the $n$ notions of
$d$--extremality, there is no need to distinguish between them.
Indeed,  for measures  $\mu$ on $\R^n$ it makes perfect sense to
redefine the standard notion of extremal  and say that {\em $\mu$
is extremal} if $\mu$ is $d$--extremal for all $d \in \{0, \ldots
n-1\}$. Recall, that the standard notion  only requires that $\mu$
is $d$--extremal for $d=0$ and $d=n-1$.

\noindent In view of  \cite{Laurent}, it is  natural to extend the
notion of $d$-extremality to the inhomogeneous setup. Here an
approximating $d$-dimensional rational plane is shifted by an
appropriately scaled transversal vector parameterised by
$\bm\theta\in \R^{n-d}$  -- the inhomogeneous part.  Leaving aside
the details, let $\omega_d(\vv x,\bm \theta)$ denote  the
inhomogeneous Diophantine exponent that arises in this way. Then,
given a measure $\mu$ on $\R^n$ we say that $\mu$ is {\em
inhomogeneously $d$--extremal} if for all $\bm \theta \in
\R^{n-d}$
$$
\text{  $\omega_d(\vv x, \bm\theta )=(d+1)/(n-d)$ \qquad for
$\mu$-almost all $\vv x \in \R^n$. }
$$
\noindent In view of the above case for redefining  the
(homogeneous)  notion of extremal for measures  $\mu$ on $\R^n$,
it would be quite natural to say that  {\em $\mu$ is
inhomogeneously extremal} if $\mu$ is inhomogeneously
$d$--extremal for all $d \in \{0, \ldots n-1\}$. The definition
given in \S\ref{inhomtheory} only requires that $\mu$ is
inhomogeneously $d$--extremal for $d=0$ and $d=n-1$. In any case,
the problem of establishing the inhomogeneous generalisation of
Theorem \ref{t10} now arises.

\noindent{\bf Conjecture 1} \   {\em Any friendly measure on
$\R^n$ is inhomogeneously $d$-extremal for all $d \in \{0, \ldots
n-1\}$.}

\medskip

\noindent Unlike the homogeneous case, the Schmidt-Laurent
transference inequalities  are  not applicable and the $n$ notions
of inhomogeneously $d$--extremality are not necessarily equivalent.
Consequently, each one needs to   be considered separately. Clearly,
Theorem \ref{t2} implies the desired statement for $d=0$ and $d=n-1$
but they do not imply the statement for  $d$ between these extreme
values as in the homogeneous case. For measures $\mu$ on $\R^n$ one
can define the notion of {\em $d$-contracting} by replacing the
planes $\cL_{\vv a,\vv b}$ appearing  in the definition of
contracting (Definition~\ref{def3}) with $d$-dimensional planes.
Conjecture 1 would then follow on establishing the following
generalization of Theorem~\ref{t1}A subject to showing that friendly
measures are $d$-contracting.

\medskip

 \noindent{\bf Conjecture 2} \   {\em Let $\mu$ be a measure on
$\R^n$ and $d \in \{0, \ldots n-1\}$. If  $\mu$ is $d$-contracting
almost everywhere then \vspace{2ex}

\centerline{$\mu$ is $d$-extremal $\iff$ $\mu$ is inhomogeneously
$d$-extremal.}

}

\medskip

\noindent In order to establish Conjecture 2, it would be natural
to follow the recipe used  to  establish Theorem \ref{t1}A. This
involves obtaining upper and lower bounds for $\omega_d(\vv
x,\bm\theta)$ separately. In short, the upper bound would follow
on applying the Inhomogeneous Transference Principle. Of course,
this is  subject to  being able to reformulate the upper bound
problem analogues to that of \S\,\ref{sec:1} and verifying the
contracting axioms of  \S\,\ref{HIT}. The lower bound would follow
on establishing an analogue of Theorem~BL or at least an analogue
of Cassels' Theorem~VI in \cite[Chapter 5]{Cassels-1957} within
the setting of approximating points in $\R^n$  by $d$-dimensional
planes. More precisly, concerning the former we seek a lower bound
for $\omega_d(\vv x,\bm\theta)$ in terms of the uniform
homogeneous exponent $\widehat\omega_{d'}(\vv x)$ for some $0 \leq
d' \leq n-1$. This is an intriguing problem in its own right and
we suspect that the following  represents the precise
relationship.

\bigskip

\noindent{\bf Conjecture 3} \   {\em Let $\vv x\in\R^n$ and
$d\in\{0,\dots,n-1\}$. Then for all $\bm \theta \in \R^{n-d}$
$$\omega_d(\vv x,\bm\theta)\ \ge\
\frac{1}{\widehat\omega_{n-1-d}(\vv x)}.$$ }

\medskip

\noindent Note that Conjecture 3 coincides with the left hand side
inequality of (\ref{e:037}) in the extreme cases $d=0$ and
$d=n-1$.

\bigskip\bigskip\noindent\textbf{Beyond  extremality in $\R^n$: system of linear forms.}
The theory of inhomogeneous extremality for measures on $\Rmn$
corresponding to a genuine system of linear forms is not covered by
Theorem~\ref{t2} or indeed Conjecture 1 above. By genuine we simply
mean that both $m$ and $n$ are strictly greater than one and so we
are outside of the dual and simultaneous theories. Naturally it
would be highly desirable to establish a general result for measures
on $\Rmn$ analogues to Theorem \ref{t2} or simply Theorem \ref{t3}.
However, such a result is currently non-existent even in the
homogeneous setting for manifolds and constitutes a key open
problem. Indeed, the homogeneous `manifold' problem first eluded to
in \cite[\S6.2]{Kleinbock-Margulis-98:MR1652916} is formally stated
in \cite[Question~35]{Gorodnik-07:MR2261070}. Part of the issue lies
in determining the right formulation of non-degeneracy and more
generally friendly condition. However, in view of Theorem~\ref{t1}
any progress on the homogeneous extremality problem can be
transferred over to the inhomogeneous setting.

\bigskip\bigskip\noindent\textbf{Beyond extremality: Khintchine-Groshev type results.}
The inhomogeneous extremality results obtained in this paper
constitute the first step towards developing a coherent
inhomogeneous theory for manifolds in line with the homogeneous
theory \cite{Beresnevich-02:MR1905790,
Beresnevich-Bernik-Kleinbock-Margulis-02:MR1944505,
Bernik-Kleinbock-Margulis-01:MR1829381}.  It would be desirable to
adapt the techniques of this paper to obtain the inhomogeneous
analogues of the convergence Khintchine-Groshev type results.

\bigskip\bigskip\noindent\textbf{Beyond Euclidean spaces\,: $\C$,
$\Q_p$ and $S$-arithmetic.} In a nutshell, Theorem~\ref{t1} enables
us to transfer homogeneous extremality statements for measures on a
Euclidean space to inhomogeneous statements. It would be desirable
to obtain analogous results within the $p$-adic, complex or more
generally the $S$-arithmetic setup. To this end, we refer the reader
to the papers
\cite{Beresnevich-Bernik-Kovalevskaya-05:MR2124042,Beresnevich-Kovalevskaya-03:MR1993537,
Ghosh-07:MR2321374,
Kleinbock-04:MR2094125,Kleinbock-Tomanov-07:MR2314053,
Mohammadi-Golsefidy-preprint} and references within for the various
homogeneous results. The Inhomogeneous Transference
Principle (Theorem \ref{t5}) is applicable within these
non-Euclidean settings and provides a natural path for obtaining the
desired inhomogeneous generalisations.

\bigskip\bigskip\noindent\textbf{Beyond rigid inhomogeneous approximation.}
Within the  context of the Diophantine approximation problems
addressed in this paper -- namely that of inhomogeneous
extremality --  the inhomogeneous part  is arbitrary but always
fixed. However, for various Diophantine approximation problems
there are often major  advantages in treating the inhomogeneous
part as a variable.  For example, in (\ref{e:001}) we may consider
the situation in which  $\bm\theta \in \R^m$ is a  function of the
approximated point $\xmn \in \Rmn $ under consideration.
Geometrically, on allowing $\bm\theta$ to depend on $\xmn$ we
perturb the underlying approximating planes $\cL_{\vv a,\vv b}$ --
see (\ref{e:011}). These `perturbed planes'   now play the role of
the approximating objects and are not necessarily planes. However,
the fact that the approximating objects are no longer planes is
completely irrelevant when it comes to applying  the Inhomogeneous
Transference Principle (Theorem \ref{t5}). In short, as long as
the intersection and contracting properties are satisfied the
Inhomogeneous Transference Principle can be used and the nature of
the approximating objects is irrelevant.

To clarify the above discussion, we describe a conjecture that can
be treated  as an inhomogeneous  problem in which the
inhomogeneous part is  a variable. Let $\psi :\N \to \R_+$  be a
monotonic function such that $\psi(r) \to 0 $ as $r \to \infty$
and for $ x \in \R $ consider the solubility of the inequality
\begin{equation}\label{e:093}
|x^n+a_{n-1}x^{n-1}+\dots+a_1x+a_0|<\psi \big( |\vv a| \big)
\end{equation}
in integer vectors  $\vv a = (a_{n-1},\dots,a_0) \in \Z^n$. Note
that when the right hand side of (\ref{e:093}) is small, the above
inequality implies that there exists an algebraic integer close to
$x$.  The following divergence statement is due to Bugeaud
\cite{Bugeaud-?2}.

\bigskip

\noindent\textbf{Theorem B } \textit{Let $\psi$ be a monotonic
function such that $\sum_{q=1}^\infty q^{n-2}\psi(q)$ diverges.
Then, for almost all $x\in\R$ the inequality $(\ref{e:093})$ has
infinitely many solutions $\vv a \in\Z^n$. }

\bigskip

\noindent Establishing the convergence counterpart to Theorem~B
represents an open problem analogous to a problem of A.\,Baker
\cite{Baker-1966-On-a-theorem-of-Sprindzuk} settled by Bernik
\cite{Bernik-1989}.

\bigskip

\noindent{\bf Conjecture 4} \   {\em Let $\psi$ be a monotonic
function such that $\sum_{q=1}^\infty q^{n-2}\psi(q)$ converges.
Then, for almost all $x\in\R$ the inequality $(\ref{e:093})$ has
at most finitely many solutions $\vv a \in\Z^n$. }

\medskip

\noindent In the special case that  $\psi(q)=q^{-n+1-\ve}$ with
$\ve>0$, the above conjecture can be viewed as the algebraic
integers analogue of Mahler's problem \cite{Mahler-1932b}. In terms
of establishing Conjecture 4, our approach is to view it as a
problem in inhomogeneous Diophantine approximation in which the
inhomogeneous part $\theta \in \R$ is a variable.  This is simple
enough to do! With reference to inequality $(\ref{e:093})$, we let
$\theta$ be $x^n$. Thus, Conjecture~4 is equivalent to a `perturbed'
inhomogeneous Diophantine approximation problem restricted to the
Veronese curve $\cV_{n-1}$ in $\R^{n-1}$. We claim that the
perturbed inhomogeneous problem can be treated via the methods
introduced in this paper.  We intend to return to Conjecture 4  in a
forthcoming paper.

\vspace{6ex}

\noindent{\em Acknowledgements. } SV would like to thank the
dynamic duo Iona and Ayesha for prolonging his youth -- of course
only in the mind!  Also, happy number six girls and ``oo ee a a a
walla walla bing bang hey!'' to you too. VB would like to thank
Vasili~Bernik for an interesting discussion that has triggered
this work. Also VB and SV wish him happy {\sc sixty one} -- many
happy returns of the day, Basil!

\bigskip
\bigskip

{\footnotesize

\begin{minipage}{\textwidth}
{\sc Victor Beresnevich \\  University of York, Heslington, York, YO10 5DD, England}\\
{\it E-mail address}\,: vb8@york.ac.uk
\end{minipage}

\bigskip

\begin{minipage}{\textwidth}
{\sc Sanju Velani \\  University of York, Heslington, York, YO10 5DD, England}\\
{\it E-mail address}\,: slv3@york.ac.uk
\end{minipage}

\bigskip
\bigskip

  \def\cprime{$'$}

}

\end{document}